\documentclass[final]{siamltex}
\usepackage{mathrsfs}
\usepackage{amsmath}
\usepackage{amsfonts}
\usepackage{amssymb}
\usepackage{graphicx}
\usepackage{geometry}
\usepackage{placeins}
\usepackage{mathrsfs}
\usepackage{ulem}
\usepackage{color}
\usepackage{verbatim}
\usepackage{subfigure}
\usepackage{graphicx}

\newtheorem{rem}{Remark}[section]

\begin{document}

\title{Novel ensemble algorithms for random two-domain parabolic problems}
\author{
Yizhong Sun \thanks{%
School of Mathematical Sciences, East China Normal University,
Shanghai. \texttt{bill950204@126.com}.}
\and Jiangshan Wang \thanks{%
School of Mathematical Sciences, East China Normal University,
Shanghai. \texttt{52215500018@stu.ecnu.edu.cn}.}
\and Haibiao Zheng \thanks{%
School of Mathematical Sciences, East China Normal University,
Shanghai Key Laboratory of Pure Mathematics and Mathematical
Practice, Key Laboratory of Advanced Theory and Application in Statistics and Data Science (East China Normal University), Shanghai, P.R. China. \texttt{hbzheng@math.ecnu.edu.cn}.
Partially supported by  the Ministry of Science and Technology (Grant No. 2022YFA1004403), Science and Technology
Commission of Shanghai Municipality (Grant Nos. 22JC1400900, 21JC1402500, 22DZ2229014) and NSF of China (Grant No. 11971174). }}

 \maketitle

\begin{abstract}
In this paper, three efficient ensemble algorithms are proposed for fast-solving the random fluid-fluid interaction model. Such a model can be simplified as coupling two heat equations with random diffusion coefficients and a friction parameter due to its complexity and uncertainty. We utilize the Monte Carlo method for the coupled model with random inputs to derive some deterministic fluid-fluid numerical models and use the ensemble idea to realize the fast computation of multiple problems. Our remarkable feature of these algorithms is employing the same coefficient matrix for multiple linear systems, significantly reducing the computational cost. By data-passing partitioned techniques, we can decouple the numerical models into two smaller sub-domain problems and achieve parallel computation. Theoretically, we derive that both algorithms are unconditionally stable and convergent. Finally, numerical experiments are conducted not only to support the theoretical results but also to validate the exclusive feature of the proposed algorithms.
\end{abstract}


\section{Introduction}
In mathematical physics, models of atmosphere-ocean interaction are generally described by two incompressible Newtonian fluids coupled with some appropriate interface conditions \cite{atmosphere-ocean1,atmosphere-ocean2,atmosphere-ocean3,atmosphere-ocean4}. In recent years, a lot of attention has been paid to the multi-domain, multi-physics coupled problems \cite{multi-domain1,multi-domain2,multi-domain3,multi-domain4,multi-domain5,multi-domain6}, especially, many novel numerical simulations are proposed for the fluid-fluid interaction problems \cite{heat-linear, heat-nolinear, fluid-fluid3, fluid-fluid4}. In \cite{heat-linear}, a simplified atmosphere-ocean interaction model with deterministic friction parameter $\kappa$ is considered as a linear heat-heat coupled system. To decouple such a multi-domain, multi-physics system
naturally and achieve parallel computation, Connors et al. \cite{heat-linear} presented the implicit-explicit partitioned method and data-passing partitioned method, which are both first-order in time, fully discrete methods.
More importantly, the highlight of \cite{heat-linear} is that the data-passing partitioned method is unconditionally stable and convergent. Then for the atmosphere-ocean interaction model with some nonlinear interface conditions, based on the existing research \cite{heat-linear}, Connors et al. \cite{heat-nolinear} further developed an unconditionally stable method by using geometric averaging for nonlinear terms.

 Due to the inaccuracy of observation data, the complexity of the atmosphere-ocean coupling, or the introduction of additional uncertainty sources, the friction parameter $\kappa$ \cite{fluid-fluid4} and diffusion coefficients $\nu_1, \nu_2$ \cite{wang zhu} are physically impossible to determine, which can only give an approximate value range or meet a certain probability distribution. Therefore, exploring the numerical simulation of such problems with uncertain inputs is necessary. Typically, we regard the uncertain parameters as random functions determined by the basic random fields with specific covariance structures (usually experimentally constructed). For simplicity, we first focus on the linear heat-heat coupled systems with three random coefficients $\kappa, \nu_1$, and $\nu_2$ as simplified fluid-fluid models. To make the analysis more transparent and understandable, we start with the simple case and then discuss the advanced case, so our analytical technique is divided into two cases:
\begin{itemize}
     \item CASE 1: Diffusion coefficients $\nu_1, \nu_2$ are deterministic parameters, we only discuss the problem in which the friction parameter $\kappa$ is random. 
    \item CASE 2: We discuss a more advanced problem that the friction parameter $\kappa$ and diffusion coefficients $\nu_1, \nu_2$ are both random. Note that $\kappa$ is a random function that is independent of space and time, but $\nu_1, \nu_2$ are the random functions of space and time.
\end{itemize}

One of the most popular approaches to solving random problems is the Monte Carlo method \cite{wang zhu, random PDEs2}, which transforms random PDEs into a series of traditional PDEs, then the existing standard numerical methods are implemented and allowed to be used. However, to better estimate the uncertainty and sensitivity of the solution, more samples are needed, which will lead to a slow convergence rate. Then, many linear equations with different stiffness matrices are formed, which inevitably requires a lot of computing costs. To overcome such computational challenges and improve its efficiency, a fast ensemble time-stepping algorithm was proposed in \cite{ensemble1}, where $J$ (the number of samples) Navier-Stokes equations with different initial conditions and forcing terms are simulated. Jiang et al. \cite{ensemble1} skillfully put forward the ensemble idea, at each time step, to solve the linear systems with one shared coefficient matrix and $J$ right hand sides. It is clear that the multiple linear systems with a shared common coefficient matrix only need to use once efficient iterative solvers, which significantly reduces the storage requirements and computational costs. For uncertainties in initial conditions and forcing terms, the ensemble algorithm has been extensively studied \cite{wang zhu, random PDEs2, ensemble1, ensemble2, ensemble3, ensemble4, ensemble5}. In \cite{wang zhu}, Wang et al. first used the ensemble algorithm to solve the parabolic problems with random parameters depending on space and time, which can be an excellent guide for our CASE2. In addition, there are many studies on more complex stochastic partial differential equations \cite{spde1, spde2,spde3}.

In this paper, we follow the ensemble idea to develop three novel algorithms for the random heat-heat coupled model. We obtained the following important results:
\begin{itemize}
    \item For the random heat-heat model, the convergence rate is rigorously demonstrated, which not only reflects the disadvantage of the Monte Carlo method with a very slow convergence rate of $1/\sqrt{J}$, but also shows the error contribution of the finite element approximation. 
     \item For CASE 1: We propose the ensemble algorithm for solving the heat-heat model only related to the random friction parameter $\kappa$. The unconditional stability and convergence of the proposed algorithm are strictly analyzed. More importantly, the random $\kappa$ does not need a small perturbation constraint, which can better simulate the real problems with large disturbance $\kappa$, usually between $10^{-3}$ and $10^{3}$.
     \item For CASE 2: Two ensemble algorithms are proposed to fast solve more advanced heat-heat models. The first algorithm can be regarded as the combination of CASE 1 and \cite{wang zhu}, but such analysis of the algorithm is more complex than the single domain problem in \cite{wang zhu}. Moreover, we do not need to strengthen any constraint conditions to get the unconditionally stable and convergent. The random diffusion coefficients $\nu_1,\nu_2$ may be the functions related to not only space but also time, so the coefficient matrix can be also changed with the time step update. To further reduce the calculation costs, we further skillfully average the coefficients $\nu_1,\nu_2$ in time on the basis of the first algorithm. No restrictions need to be added or strengthened for the unconditional stability and convergence of the second algorithm.
\end{itemize}

 The rest of this paper is organized as follows. In Section 2, we introduce the random heat-heat coupled model. Also, some notations and mathematical preliminaries are described. In Section 3, we propose an efficient ensemble algorithm for CASE 1 and prove its unconditional stability and convergence. Two ensemble algorithms are proposed for CASE 2, and then the stability and error estimation of both algorithms are demonstrated in Section 4. Finally, numerical tests are presented to illustrate the effectiveness of the proposed scheme in Section 5.

\section{Notations and preliminaries}
In this paper, we will consider two heat equations coupled by some interface conditions that allow energy to transform back and forth across one interface. Let the domain $\Omega \subset \mathbb{R}^d \hspace{0.5mm} (d=2,3)$ have convex, polygonal subdomains $\Omega_i \hspace{0.5mm} (i=1,2)$ with the interface $\Gamma=\partial{\Omega}_1\cap\partial{\Omega}_2=\bar{\Omega}_1\cap\bar{\Omega}_2$ and outer boundary of each subdomain $\Gamma_i=\partial\Omega_i~\backslash~\Gamma$. Denote by $\hat{n}_1$ and $\hat{n}_2$ the unit outward normal vectors on $\partial\Omega_1$ and $\partial\Omega_2$, respectively. Fig. 2.1 is a simple 2D diagram to help understand.
\begin{figure}[htbp]\label{domain0}
\centering
\includegraphics[width=90mm,height=35mm]{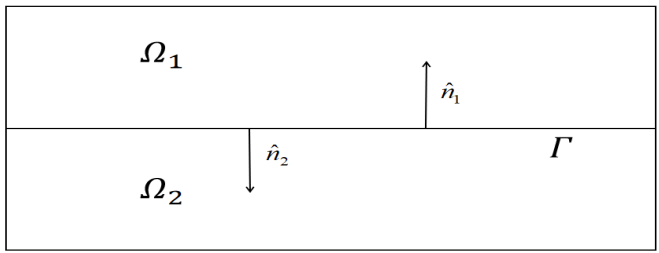}
\caption{Example sub-domains, adjoined by an interface $\Gamma$. }
\end{figure}

The two heat equations coupled by the jump interface conditions are shown as follows:
\begin{eqnarray}
u_{i,t}-\nu_i\Delta u_i&=&f_i  \ \ \ \ \ \ \qquad \ \ \ \mathrm{in}~\Omega _{i}\times[0,T],  \label{heat1}\\
-\nu_i\nabla u_i \cdot\hat{n}_i&=&\kappa(u_i-u_k) \ \ \ \ \mathrm{on}~\Gamma\times[0,T], \ i,k=1,2, \ i\neq k, \label{heat2} \\
u_i({\bf x},0)&=&u^0_i({\bf x})    \ \ \ \quad \ \  \ \ \mathrm{in}~\Omega _{i}\times[0,T],   \label{heat3}    \\
u_i({\bf x}, t)&=&0 \ \ \qquad \qquad \    \ \mathrm{on}~\Gamma_{i}=\partial\Omega_i~\backslash~\Gamma,\label{heat4}
\end{eqnarray}
where $\nu_i \in L^2(W^{1,\infty}(\Omega);0, T)$ denotes the diffusion coefficient, $f_i \in L^2(H^{-1}(\Omega);0, T)$ is the sink/source term, and $\kappa>0$ means the friction parameter, which is calculated in practice from bulk flux formulae \cite{fluid-fluid4}. Moreover, the initial condition $u_i^0$ need to satisfy $u_i^0 \in H_0^1(\Omega)\cap H^{p+1}(\Omega)$ with $p\ge 1$.

Due to the inaccuracy of observation data and the complexity of the atmosphere-ocean coupling, the diffusion coefficients $\nu_1, \nu_2$, and the friction parameter $\kappa$ are physically impossible to determine. Moreover, some additional uncertainty sources may need to be introduced to get $\kappa$ \cite{fluid-fluid4}. This forces us to further study the random heat-heat coupled model with three random coefficients $\kappa, \nu_1$, and $\nu_2$. Let $(\Pi, \mathcal{F}, \mathcal{P})$ be a complete probability space, where $\Pi$ is the set of outcomes, $\mathcal{F} \subset 2^{\Pi}$ is the $\sigma$-algebra of events, and $\mathcal{P}: \mathcal{F}\to [0,1]$ is a probability measure. That is, the random heat-heat coupled system reads: find the random function $u_i$: $\Pi \times \Omega_i \times [0, T] \to \mathbb{R}$ satisfying $\mathcal{P} - a.e.$
\begin{eqnarray}
u_{i,t}-\nu_i(\omega, {\bf x}, t)\Delta u_i&=&f_i(\omega, {\bf x}, t) \ \  \qquad \ \ \ \ \mathrm{in}~\Pi \times \Omega_i \times [0,T],  \\
-\nu_i(\omega, {\bf x}, t)\nabla u_i \cdot\hat{n}_i&=&\kappa(\omega)(u_i-u_k) \  \ \ \mathrm{on}~\Pi \times \Gamma \times [0,T], \hspace{0.5mm} i,k=1,2, \hspace{0.5mm} i\neq k,  \\
u_i(\omega, {\bf x},0)&=&u^0_i(\omega, {\bf x})   \ \ \ \  \ \qquad \  \ \ \mathrm{in}~\Pi \times \Omega_i \times[0,T],      \\
u_i(\omega, {\bf x}, t)&=&0   \ \ \qquad \qquad  \   \qquad \  \ \mathrm{on}~\Pi \times \Gamma_{i}=\partial\Omega_i~\backslash~\Gamma,
\end{eqnarray}
where the diffusion coefficient $\nu_i$ and source force $f_i$ hold: $\Pi \times \Omega_i \times [0, T] \to \mathbb{R}$, and $u_i^0: \Pi \times \Omega_i \to \mathbb{R}$, which is assumed to have continuous and bounded covariance functions. 

The Monte Carlo method is one of the most classical approaches to solving random PDEs. The main idea of this method is to obtain the identically distributed approximations of the solutions by repeated sampling of input parameters, and to use standard numerical methods for the corresponding deterministic PDEs. Then the approximate solutions are further analyzed to yield useful statistical information. The computation procedure consists of the following step:

1. Choose a set of independently, identically distributed (i.i.d) samples for the random friction parameter $\kappa(\omega_j)$ and diffusion coefficients $\nu_i(\omega_j, \cdot, \cdot)$, $j=1,\cdots, J$;

2. Apply standard numerical method to solve for approximate solutions $u_{1,j}^{n+1}$ and $u_{2,j}^{n+1}$, $j=1,\cdots, J$;

3. Approximate expectation $E[{\bf u}]$ by averaging the statistical information of the output: $E[{\bf u}]\approx \frac{1}{J}\sum\limits_{j=1}^J{\bf u}(\omega_j,\cdot,\cdot)$.

\begin{rem}
When solving the corresponding deterministic PDEs with a standard numerical method, we need to solve the following linear systems
\begin{eqnarray*}
A_{1,j}\left[\begin{array}{l l l}
{u}_1(\omega_j, {\bf x}, t) \\
\end{array}\right]=\left[RHS_{1,j}({\bf x})\right], \quad
A_{2,j}\left[\begin{array}{l l l}
{u}_2(\omega_j, {\bf x}, t) \\
\end{array}\right]=\left[RHS_{2,j}({\bf x}) \right].
\end{eqnarray*}
As shown above, to better estimate the uncertainty and sensitivity of the solutions, more samples are needed. Then a large number of linear equations with different stiffness matrices are formed, which will increase the computational cost.
\end{rem}

Suppose $\bar\nu_i=\frac{1}{J}\sum\limits_{j=1}^J\nu_{i}(\omega_j,{\bf x},t)$, then the following two conditions are valid:

(\romannumeral1) There exists a positive constant $\theta$ such that, for any $t \in [0,T]$,
\begin{eqnarray}
P\left\{ \omega \in \Omega; \min\limits_{{\bf x}\in \bar{\Omega}} \bar\nu_i(\omega, {\bf x},t)\ge \theta\right\}=1.
\end{eqnarray}

(\romannumeral2) There exist positive constants $\theta_-$ and $\theta_+$ such that, for any $t\in [0,T]$,
\begin{eqnarray}
P\left\{\omega_j \in \Omega; \theta_-\le |\nu_{i}(\omega_j, {\bf x},t)-\bar{\nu}_i({\bf x},t)|_{\infty}\le\theta_+ \right\}=1.
\end{eqnarray}
\begin{theorem}
Under conditions (\romannumeral1) and (\romannumeral2), suppose $\theta>\theta_+$ and $ \Phi^n=\frac{1}{J}\sum\limits_{j=1}^J{\bf u}_j^n$, then we have the error estimation of (2.5)-(2.8) as follows
\begin{eqnarray}
&&E[\|E[{\bf u}^{N+1}]-\Phi^{N+1}\|^2]+\kappa_{max}\Delta t E[\|E[{\bf u}^{N+1}]-\Phi^{N+1}\|_\Gamma^2]\nonumber\\
&&\hspace{1mm}+\theta_-\Delta t E[\|\nabla E[{\bf u}^{N+1}]-\nabla\Phi^{N+1}\|^2]+(\theta-\theta_+)\Delta t\sum\limits_{n=0}^N E[\|\nabla E[{\bf u}^{n+1}]-\nabla\Phi^{n+1}\|^2]\nonumber\\
&&\le \frac{1}{J}\left(E[\|{\bf u}_j^0\|^2]+C\Delta tE[\|\nabla {\bf u}_j^0\|^2]+C\Delta t E[\|{\bf u}_j^0\|_\Gamma^2]+C\Delta t\sum\limits_{n=0}^N E[\|{\bf f}_j^{n+1}\|_{-1}^2]\right) \nonumber\\
&&\hspace{3mm}+C(\Delta t^2+h^{2l}).
\end{eqnarray}
\end{theorem}

\begin{proof}
The proof of the error estimate conclusion can be divided into two parts. One is the standard error estimate of the Monte Carlo method, and the analysis process can refer to \cite{wang zhu}. The other part of the error estimate comes from the ensemble algorithm, and its proof will be given later (Theorem 4.2).
\end{proof}

When regarding random friction parameter $\kappa(\omega_j)$ as $\kappa_j$, and random diffusion coefficients $\nu_i(\omega_j, {\bf x}, t)$ as $\nu_{i,j}$, similarly denoting sink/source term $f_i(\omega_j,{\bf x},t)$ as $f_{i,j}$, we have an ensemble of $J$ heat-heat coupled systems corresponding to $J$ different parameters sets $(f_{i,j}, \kappa_j,\nu_{i,j}), \hspace{0.5mm} i=1,2, \hspace{0.5mm} j=1,\cdots, J$ as follows:
\begin{eqnarray}
 u_{i,j,t}-\nu_{i,j}\Delta u_{i,j}&=&f_{i,j} \ \ \ \ \qquad \qquad \ \mathrm{in}~\Omega _{i} \times [0,T], \label{random1.1}\\
-\nu_{i,j}\nabla u_{i,j} \cdot\hat{n}_i&=&\kappa_j(u_{i,j}-u_{k,j}) \ \ \  \mathrm{on}~\Gamma \times [0,T], i,k=1,2, \hspace{0.5mm} i\neq k, \label{random1.2} \\
u_{i,j}({\bf x},0)&=&u^0_{i,j}({\bf x})   \ \ \qquad  \ \ \ \ \ \mathrm{in}~\Omega _{i} \times[0,T],    \label{random1.3}    \\
u_{i,j}({\bf x}, t)&=&0   \ \ \ \qquad \qquad \quad \  \ \mathrm{on}~\Gamma_{i}=\partial\Omega_i~\backslash~\Gamma.\label{random1.4}
 \end{eqnarray}

Next, we introduce some notations, Sobolev spaces, and norms. Denote the $L^2(\Omega)$ norm and inner product by $\|\cdot\|$ and $(\cdot,\cdot)$, respectively, and $L^2(\Gamma)$ norm on the interface by $\|\cdot\|_{\Gamma}$. Let $X_i:=\{v_i\in H^1(\Omega_i): v_i=0~\mathrm{on}~\Gamma_i$\} ($i=1,2$). 
 We define $X=X_1\times X_2=\{{\bf v}=(v_1,v_2): v_i \in X_i, \hspace{0.5mm} i=1,2\}$ and $L^2(\Omega)=L^2(\Omega_1)\times L^2(\Omega_2)$ as two product spaces for the global domain.
 Moreover, the product spaces $L^2(\Omega)$ and $X$ are equipped with the following norms: $\|{\bf v}\|=\sum_{i=1,2}(\int_{\Omega_i}|v_i|^2 dx)^{1/2},$ $\|{\bf v}\|_X=\sum_{i=1,2}(\int_{\Omega_i}(|v_i|^2+|\nabla v_i|^2)dx)^{1/2}.$ The space $H^{-1}(\Omega)$ denotes the dual space of bounded linear functions defined on $H^1_0(\Omega)$, a norm for $H^{-1}(\Omega)$ is given by
$\|{\bf f}\|_{-1}=\sup\limits_{0\neq{\bf v}\in H_0^1(\Omega)} \frac{({\bf f},{\bf v})}{\|\nabla {\bf v}\|}$, where ${\bf f}=(f_1,f_2)$.

With the above notation, the weak formulation for (\ref{random1.1})-(\ref{random1.4}) can be written as follows: for $i,k=1,2, \hspace{0.5mm} i\neq k$, find $u_{i,j}:[0,T]\to X_i $ satisfying
\begin{eqnarray}
(u_{i,j,t},v_i)+(\nu_{i,j}\nabla u_{i,j},\nabla v_i)+\int_\Gamma {\kappa_j}(u_{i,j}-u_{k,j})v_i\mathrm{d}s=(f_{i,j},v_i) \quad \forall v_i\in X_i, \label{week1}
\end{eqnarray}

The natural monolithic weak formulation for (\ref{random1.1})-(\ref{random1.4}), found by summing (\ref{week1}) over $i=1,2$, is shown as below: find ${\bf u}_j:[0,T] \to X$ satisfying
\begin{eqnarray}\label{Moweakform}
({\bf u}_{j,t}, {\bf v})+({\nu}_j \nabla{\bf u}_j, \nabla{\bf v})+\int_\Gamma\kappa_j[{\bf u}_j][{\bf v}]\mathrm{d} s=({\bf f}_j, {\bf v}), \quad \forall {\bf v}\in X,
\end{eqnarray}
where $[\cdot]$ denotes the jump of the indicated quantity across the interface $\Gamma$ and
\begin{eqnarray*}
({\bf u}_{j,t}, {\bf v})=\sum\limits_{i=1}^{2}({u}_{i,j,t}, {v_i}),\quad
({ \nu}_j \nabla{\bf u}_j, \nabla{\bf v})=\sum\limits_{i=1}^{2}(\nu_{i,j}\nabla{u_{i,j}}, \nabla{v_i}),\quad
({\bf f}_j, {\bf v})=\sum\limits_{i=1}^{2}({f_{i,j}}, {v_i}).
\end{eqnarray*}

Let $\mathcal{T}_i$ be a quasi-uniform triangulation of the domain $\Omega_i$ and $\mathcal{T}_h=\mathcal{T}_1\cup\mathcal{T}_2$. Define mesh size $h$ is the largest diameter of a simplex in $\mathcal{T}_i$. Take $X_{i,h}\subset X_i$ be the finite element spaces for $i=1,2$, and define $X_h=X_{1,h}\times X_{2,h} \subset X$. Divide the simulation time $T$ into $N$ smaller time intervals with $[0,T]=\mathop{\cup}\limits_{n=0}^{N-1}[t^n,t^{n+1}]$, where $t^n=n\Delta t, \Delta t=T/N$. For $j=1,\cdots, J$ and $n=0,1,2,\cdots,N$, ${\bf u}_j^n$ denotes the discrete approximation to ${\bf u}_j(t^n) \in X_h$.

Furthermore, we recall the Poincar$\acute{e}$ inequality and trace inequality. There exist positive constants $C_p$ and $C_t$ which depend only on the domain $\Omega$, such that
\begin{eqnarray*}
\|{\bf v}\|\le C_p\|\nabla {\bf v}\|, \qquad   \|{\bf v}\|_\Gamma\le C_t\|{\bf v}\|^{\frac{1}{2}}\|\nabla {\bf v}\|^{\frac{1}{2}}.
\end{eqnarray*}

The following optimal approximation properties of piecewise continuous polynomials on the quasi-uniform mesh of local $l$ are assumed for some mesh-independent constant $C$, for $i=1,2$
\begin{eqnarray*}
\inf\limits_{v\in X_{i,h}}\|u-v\|\le Ch^{l+1}|u|_{H^{l+1}},\quad
\inf\limits_{v\in X_{i,h}}\|u-v\|_{X_i}\le Ch^{l}|u|_{H^{l+1}}.
\end{eqnarray*}

In this report, we mainly commit to presenting novel numerical algorithms for the second procedure of the above Monte Carlo method, namely, the $J$ heat-heat coupled models with three random coefficients. Next, we will discuss the problem from simple to difficult in the following two sections.

\section{The ensemble algorithm for CASE 1}

We first consider CASE 1 with random friction parameter $\kappa$ only. The ensemble algorithm for CASE 1 can be proposed as follows:

\textbf{{{Ensemble Algorithm 1 (A1)}}}

Let $\Delta t>0$, $\kappa_{\max}=\max\limits_j \kappa_j$, ${\bf u}_j^0\in X_h$ and $ {\bf f}_{j}\in L^2(H^{-1}(\Omega); 0,T))$. Given ${\bf u}_j^n\in X_h$, find ${\bf u}_j^{n+1}\in X_h$ for $n=0,1,2,\cdots,N-1$ and $j=1,\ldots,J$,
\begin{eqnarray}
(\frac{u_{1,j}^{n+1}-u_{1,j}^{n}}{\Delta t},v_1)+(\nu_1\nabla u_{1,j}^{n+1},\nabla v_1)+\int_\Gamma {\kappa_{\max}}(u_{1,j}^{n+1}-u_{2,j}^n) v_1\mathrm{d}s \nonumber\\
+\int_\Gamma (\kappa_j-{\kappa_{\max}})(u_{1,j}^n-u_{2,j}^n) v_1 \mathrm{d}s=(f_{1,j}^{n+1},v_1) \ \ \ \ \forall v_1\in X_{1,h},\label{alrorithm1.1}\\
(\frac{u_{2,j}^{n+1}-u_{2,j}^{n}}{\Delta t},v_2)+(\nu_2\nabla u_{2,j}^{n+1},\nabla v_2)+\int_\Gamma {\kappa_{\max}}(u_{2,j}^{n+1}-u_{1,j}^n) v_2\mathrm{d}s\nonumber\\
+\int_\Gamma (\kappa_j-{\kappa_{\max}})(u_{2,j}^n-u_{1,j}^n) v_2 \mathrm{d}s=(f_{2,j}^{n+1},v_2) \ \ \ \  \forall v_2\in X_{2,h}.\label{alrorithm1.2}
\end{eqnarray}


\begin{rem}
At each time step, we only compute the solutions of two linear systems, which share the same coefficient matrix:
\begin{eqnarray*}
A_1\left[\begin{array}{l l l}
{u}_{1,1} \\
\end{array}
|\cdots|
\begin{array}{l}
{u}_{1,J} \\
\end{array}\right]&=&\left[RHS_{1,1}|\cdots| RHS_{1,J}\right], \\
A_2\left[\begin{array}{l l l}
{u}_{2,1} \\
\end{array}
|\cdots|
\begin{array}{l}
{u}_{2,J} \\
\end{array}\right]&=&\left[RHS_{2,1}|\cdots| RHS_{2,J}\right].
\end{eqnarray*}
Hence the coefficient matrices $A_1$ and $A_2$ only need to use once efficient iterative solves or direct solvers such as $LU$ factorization for fast computation, which can reduce storage and computation time. Moreover, the proposed algorithm can decouple the original problem into two independent sub-physical problems and will be solved in parallel.
\end{rem}

Next, supposing $\nu_{\min}:=\min_{{\bf x}\in \bar{\Omega}} \{ \nu_1({\bf x},t), \nu_2({\bf x}, t) \}$, we can establish the stability of the approximations in Ensemble Algorithm 1 (A1) here.
\begin{theorem}(A1 Stability)
Suppose that ${\bf f}_j\in L^2(H^{-1}(\Omega);0,T)$ and ${\bf u}_j^{n+1}\in X_h$ satisfy (\ref{alrorithm1.1})-(\ref{alrorithm1.2}) for each $n\in \{0,1,2,\cdots,N-1\}$, then there exists  a generic positive constant $C$ independent of $h, \Delta t, J$, such that the numerical solution ${\bf u}_j^{n+1}$ satisfies:
\begin{eqnarray}
\|{\bf u}_j^{N}\|^2&+&\nu_{\min}\Delta t\sum\limits_{n=0}^{N-1}\|\nabla{\bf u}_{j}^{n+1}\|^2+\kappa_{\max}\Delta t\|{\bf u}_j^{N}\|^2_\Gamma\nonumber\\
&&\le \|{\bf u}_j^0\|^2+C\Delta t\|{\bf u}_j^0\|_\Gamma^2+C\Delta t\sum\limits_{n=0}^{N-1}\|{\bf f}_{j}^{n+1}\|_{-1}^2.\label{stability}
\end{eqnarray}
\end{theorem}

\begin{proof}
Taking ${\bf v}={\bf u}_j^{n+1}$ in (\ref{alrorithm1.1})-(\ref{alrorithm1.2}), we have
\begin{eqnarray}
&&\left(\frac{{\bf u}_j^{n+1}-{\bf u}_j^n}{\Delta t}, {\bf u}_j^{n+1} \right)+(\nu \nabla{\bf u}_{j}^{n+1},\nabla{\bf u}_{j}^{n+1})+\kappa_{\max}\|{\bf u}_j^{n+1}\|_\Gamma^2+\int_\Gamma(\kappa_j-\kappa_{\max})u_{1,j}^n u_{1,j}^{n+1}\mathrm{d}s\nonumber\\
&&\hspace{2mm}
+\int_\Gamma(\kappa_j-\kappa_{\max})u_{2,j}^n u_{2,j}^{n+1}\mathrm{d}s
-\int_\Gamma\kappa_j u_{2,j}^n u_{1,j}^{n+1}\mathrm{d}s-\int_\Gamma\kappa_j u_{1,j}^n u_{2,j}^{n+1}\mathrm{d}s=({\bf f}_{j}^{n+1},{\bf u}_{j}^{n+1}).\label{stability1}
\end{eqnarray}
Using the Young's inequality and $2(a-b, a)=a^2+(a-b)^2-b^2$, we get
\begin{eqnarray}
&&\frac{1}{2\Delta t}(\|{\bf u}_j^{n+1}\|^2-\|{\bf u}_j^n\|^2)+{\nu_{\min}}\|\nabla{\bf u}_{j}^{n+1}\|^2+\kappa_{max}\|{\bf u}_j^{n+1}\|_\Gamma^2\nonumber\\
&&\hspace{8mm}\le\int_\Gamma(\kappa_{\max}-\kappa_j)u_{1,j}^n u_{1,j}^{n+1}\mathrm{d}s+\int_\Gamma(\kappa_{\max}-\kappa_j)u_{2,j}^n u_{2,j}^{n+1}\mathrm{d}s  \nonumber\\
&&\hspace{12mm}+\int_\Gamma\kappa_j u_{2,j}^n u_{1,j}^{n+1}\mathrm{d}s+\int_\Gamma\kappa_j u_{1,j}^n u_{2,j}^{n+1}\mathrm{d}s+({\bf f}_{j}^{n+1},{\bf u}_{j}^{n+1}).\label{stability2}
\end{eqnarray}
Applying the Cauchy-Schwarz and Young's inequalities on the right hand side (RHS), we have
\begin{eqnarray}
&&\int_\Gamma(\kappa_{\max}-\kappa_j)u_{1,j}^n u_{1,j}^{n+1}\mathrm{d}s+\int_\Gamma(\kappa_{\max}-\kappa_j)u_{2,j}^n u_{2,j}^{n+1}\mathrm{d}s
+\int_\Gamma\kappa_j u_{2,j}^n u_{1,j}^{n+1}\mathrm{d}s+\int_\Gamma\kappa_j u_{1,j}^n u_{2,j}^{n+1}\mathrm{d}s\nonumber\\
&&\hspace{2mm}\le\frac{\kappa_{\max}-\kappa_j}{2}\|u_{1,j}^n\|_\Gamma^2+\frac{\kappa_{\max}-\kappa_j}{2}\|u_{1,j}^{n+1}\|_\Gamma^2+\frac{\kappa_{\max}-\kappa_j}{2}\|u_{2,j}^n\|_\Gamma^2
+\frac{\kappa_{\max}-\kappa_j}{2}\|u_{2,j}^{n+1}\|_\Gamma^2\nonumber\\
&&\hspace{6mm}+\frac{\kappa_j}{2}\|u_{2,j}^n\|_\Gamma^2+\frac{\kappa_j}{2}\|u_{1,j}^{n+1}\|_\Gamma^2+\frac{\kappa_j}{2}\|u_{1,j}^n\|_\Gamma^2
+\frac{\kappa_j}{2}\|u_{2,j}^{n+1}\|_\Gamma^2\nonumber\\
&&\hspace{2mm}\le\frac{\kappa_{\max}}{2}\|{\bf u}_j^n\|_\Gamma^2+\frac{\kappa_{\max}}{2}\|{\bf u}_j^{n+1}\|_\Gamma^2,\label{stability3}\\
&&({\bf f}_{j}^{n+1},{\bf u}_{j}^{n+1})\le \|{\bf f}_{j}^{n+1}\|_{-1}\|\nabla{\bf u}_{j}^{n+1}\| \le \frac{1}{2{\nu_{\min}}}\|{\bf f}_{j}^{n+1}\|_{-1}^2+\frac{\nu_{\min}}{2}\|\nabla{ \bf u}_{j}^{n+1}\|^2.\label{stability4}
\end{eqnarray}
Substituting (\ref{stability3})-(\ref{stability4}) into (\ref{stability2}) yields
\begin{eqnarray}
\frac{1}{2\Delta t}(\|{\bf u}_j^{n+1}\|^2-\|{\bf u}_j^n\|^2)+\frac{\nu_{\min}}{2}\|\nabla{\bf u}_{j}^{n+1}\|^2+\frac{\kappa_{\max}}{2}(\|{\bf u}_j^{n+1}\|_\Gamma^2-\|{\bf u}_j^{n}\|_\Gamma^2) 
\le \frac{1}{2\nu_{\min}}\|{\bf f}_{j}^{n+1}\|_{-1}^2.\label{stability5}
\end{eqnarray}
Multiplying both sides by $2\Delta t$, and summing over $n=0,1,\cdots, N-1$, we can finally yield the unconditional stability (\ref{stability}) of A1.
\end{proof}

Next, supposing $\nu_{\max}:=\max_{\bf{x}\in \bar{\Omega}} \{ \nu_1(\bf{x},t), \nu_2(\bf{x}, t) \}$, we can estimate the approximation error of A1.
\begin{theorem}(A1 Error Estimate)
Let ${\bf u}_j(t^{n+1})$ and ${\bf u}_j^{n+1}$ be the solutions of the natural monolithic weak formulation (\ref{Moweakform}) and the algorithm A1 at time $t^{n+1}$, respectively. Assume ${\bf u}_t\in L^2(X; 0, T), {\bf u}_{tt}\in L^2(L^2(\Omega);0,T)$, for all $t\in(0.T)$. Then, for any $n\in \{0,1,2,\cdots, N-1\}$, there exists a generic positive constant $C$ independent of $h, \Delta t, J$, such that
\begin{eqnarray}\label{A1convergence}
&&\|{\bf u}_j(t^{N})-{\bf u}_j^{N}\|^2+\kappa_{\max}\Delta t\|{\bf u}_j(t^{N})-{\bf u}_j^{N}\|_\Gamma^2+{\nu}_{\min}\Delta t\sum\limits_{n=0}^{N-1}\|\nabla{\bf u}_j(t^{n+1})-\nabla{\bf u}_j^{n+1}\|^2\nonumber\\
&&\le C\bigg\{\|{\bf u}_j(0)-{\bf u}_j^0\|^2+\kappa_{\max}\Delta t\|{\bf u}_j(0)-{\bf u}_j^0\|_\Gamma^2+{\Delta t}^2\|\nabla{\bf u}_{j,t}\|_{L^2(0,T;L^2{\Omega})}^2\bigg.\nonumber\\
&&\hspace{3mm}\bigg.+{\Delta t}^2\|{\bf u}_{j,tt}\|_{L^2(0,T;L^2(\Omega))}^2+\inf\limits_{{\bf v}_j^0\in X_h}\left\{\|{\bf u}_j(0)-{\bf v}_j^0\|^2+\kappa_{\max}\Delta t\|\nabla({\bf u}_j(0)-{\bf v}_j^0)\|_\Gamma^2\right\} \bigg.\nonumber\\
&&\hspace{3mm}\bigg.+\inf\limits_{{\bf v}_j\in X_h}\|({\bf u}_j-{\bf v}_j)_t\|_{L^2(0,T;L^2(\Omega))}^2+T \max\limits_{n=0,1,\cdots,N}\inf\limits_{{\bf v}_j^n\in X_h}\|\nabla{\bf u}_j(t^n)-\nabla{\bf v}_j^n\|^2 \bigg\}.
\end{eqnarray}
\end{theorem}
\begin{proof}
Restricting test function $\bf{v}$ to $X_h$, subtracting (\ref{alrorithm1.1})-(\ref{alrorithm1.2}) from (\ref{Moweakform}), we get the following error equation
\begin{eqnarray}
&&\left({\bf u}_{j,t}(t^{n+1})-\frac{{\bf u}_j^{n+1}-{\bf u}_j^{n}}{\Delta t}, {\bf v}\right)
+\left({\nu} (\nabla{\bf u}_{j}(t^{n+1})
-\nabla{\bf u}_{j}^{n+1} ),\nabla{\bf v} \right)\nonumber\\
&&\hspace{2mm}+\int_\Gamma\kappa_j(u_{1,j}(t^{n+1})-u_{2,j}(t^{n+1})){v_1} \mathrm{d}s
+\int_\Gamma\kappa_j(u_{2,j}(t^{n+1})-u_{1,j}(t^{n+1})){v_2} \mathrm{d}s\nonumber\\
&&\hspace{2mm}-\int_\Gamma \kappa_{\max}(u_{1,j}^{n+1}-u_{2,j}^n) v_1 \mathrm{d}s
-\int_\Gamma \kappa_{\max}(u_{2,j}^{n+1}-u_{1,j}^n) v_2 \mathrm{d}s \nonumber\\
&&\hspace{2mm}-\int_\Gamma (\kappa_j-\kappa_{\max})(u_{1,j}^n-u_{2,j}^n) v_1 \mathrm{d}s
-\int_\Gamma (\kappa_j-\kappa_{\max})(u_{2,j}^n-u_{1,j}^n)v_2 \mathrm{d}s=0.\label{convergence1}
\end{eqnarray}
Define ${\bf r}_j^{n+1}={\bf u}_{j,t}(t^{n+1})-\frac{{\bf u}_j(t^{n+1})-{\bf u}_j(t^{n})}{\Delta t}$, and rearrange above terms to arrive
\begin{eqnarray}
&&({\bf r}_j^{n+1},{\bf v})+\left(\frac{{\bf u}_j(t^{n+1})-{\bf u}_j^{n+1}}{\Delta t}-\frac{{\bf u}_j(t^{n})-{\bf u}_j^{n}}{\Delta t},{\bf v}\right)+({\nu}(\nabla{\bf u}_{j}(t^{n+1})-\nabla{\bf u}_{j}^{n+1}),\nabla{\bf v})\nonumber\\
&&\hspace{2mm}+\int_\Gamma\kappa_j(u_{1,j}(t^{n+1})-u_{2,j}(t^{n+1})){v_1} \mathrm{d}s+\int_\Gamma\kappa_j(u_{2,j}(t^{n+1})-u_{1,j}(t^{n+1})){v_2}\mathrm{d}s \nonumber\\
&&\hspace{2mm}-\int_\Gamma {\kappa_{\max}}(u_{1,j}^{n+1}-u_{2,j}^n) v_1 \mathrm{d}s
-\int_\Gamma {\kappa_{\max}}(u_{2,j}^{n+1}-u_{1,j}^n) v_2\mathrm{d}s\nonumber\\
&&\hspace{2mm}-\int_\Gamma (\kappa_j-{\kappa_{\max}})(u_{1,j}^n-u_{2,j}^n) v_1 \mathrm{d}s
-\int_\Gamma (\kappa_j-{\kappa_{\max}})(u_{2,j}^n-u_{1,j}^n) v_2 \mathrm{d}s=0.\label{convergence2}
\end{eqnarray}
Some error functions are decomposed as follows
\begin{eqnarray}\label{errorfunc}
{\bf u}_j(t^n)-{\bf u}_j^n=({\bf u}_j(t^n)-{\bf v}_j^n)+({\bf v}_j^n-{\bf u}_j^n)={{\eta}}_j^n+{ {\phi}}_j^n, \quad \quad \forall {\bf v}_j^n\in X_h.
\end{eqnarray}
Then, the error equation (\ref{convergence2}) can be rewritten as follows
\begin{eqnarray}
&&\frac{1}{\Delta t}({ \phi}_j^{n+1}-{ \phi}_j^n, {\bf v})+({\nu}\nabla{ \phi}_j^{n+1},\nabla{\bf v})
+\int_\Gamma\kappa_j (u_{1,j}(t^{n+1})-u_{2,j}(t^{n+1})){v_1}\mathrm{d}s\nonumber\\
&&\hspace{10mm}+\int_\Gamma\kappa_j(u_{2,j}(t^{n+1})-u_{1,j}(t^{n+1})){v_2}\mathrm{d}s
\nonumber\\
&&\hspace{10mm}-\int_\Gamma {\kappa_{\max}}(u_{1,j}^{n+1}-u_{2,j}^n) v_1\mathrm{d}s -\int_\Gamma {\kappa_{\max}}(u_{2,j}^{n+1}-u_{1,j}^n) v_2\mathrm{d}s
\nonumber\\
&&\hspace{10mm}-\int_\Gamma (\kappa_j-{\kappa_{\max}})(u_{1,j}^n-u_{2,j}^n) v_1 \mathrm{d}s-\int_\Gamma (\kappa_j-{\kappa_{\max}})(u_{2,j}^n-u_{1,j}^n) v_2 \mathrm{d}s\nonumber\\
&&\hspace{4mm}=-\frac{1}{\Delta t}({ \eta}_j^{n+1}-{ \eta}_j^n,{\bf v})-({\bf r}_j^{n+1},{\bf v})-({\nu}\nabla{ \eta}_j^{n+1},{\nabla{\bf v}}).
\label{convergence3}
\end{eqnarray}
Choosing ${\bf v}={ \phi}_j^{n+1}$ and using $2(a-b,a)=a^2+(a-b)^2-b^2$, we get
\begin{eqnarray}
&&\frac{1}{2\Delta t}(\|{ \phi}_j^{n+1}\|^2-\|{ \phi}_j^n\|^2)+{\nu}_{\min}\|\nabla {\phi}_j^{n+1}\|^2
+\int_\Gamma\kappa_j(u_{1,j}(t^{n+1})-u_{2,j}(t^{n+1})){\phi}_{1,j}^{n+1}\mathrm{d}s\nonumber\\
&&\hspace{10mm}+\int_\Gamma\kappa_j(u_{2,j}(t^{n+1})-u_{1,j}(t^{n+1})){\phi}_{2,j}^{n+1}\mathrm{d}s
\nonumber\\
&&\hspace{10mm}-\int_\Gamma {\kappa_{\max}}(u_{1,j}^{n+1}-u_{2,j}^n){\phi}_{1,j}^{n+1}\mathrm{d}s-\int_\Gamma {\kappa_{\max}}(u_{2,j}^{n+1}-u_{1,j}^n) {\phi}_{2,j}^{n+1}\mathrm{d}s\nonumber\\
&&\hspace{10mm}-\int_\Gamma (\kappa_j-{\kappa_{\max}})(u_{1,j}^n-u_{2,j}^n) {\phi}_{1,j}^{n+1} \mathrm{d}s -\int_\Gamma (\kappa_j-{\kappa_{\max}})(u_{2,j}^n-u_{1,j}^n) {\phi}_{2,j}^{n+1} \mathrm{d}s\nonumber\\
&&\hspace{4mm}\le-\frac{1}{\Delta t}({\bf \eta}_j^{n+1}-{\bf \eta}_j^n,{\bf \phi}_j^{n+1})-({\bf r}_j^{n+1},{\bf \phi}_j^{n+1})-({\nu} \nabla{\bf \eta}_j^{n+1},\nabla{\bf \phi}_j^{n+1}).\label{convergence4}
\end{eqnarray}
We will treat the interface terms of (\ref{convergence4}) in a useful way. We can add and subtract the following four terms $\int_\Gamma \kappa_{\max}(u_1(t^{n+1})-u_2(t^n)) \phi_{1,j}^{n+1} \mathrm{d}s$, $\int_\Gamma\kappa_{\max}(u_2(t^{n+1})-u_1(t^n))\phi_{2,j}^{n+1} \mathrm{d}s$, $\int_\Gamma(\kappa_j-\kappa_{\max})(u_1(t^n)-u_2(t^n))\phi_{1,j}^{n+1} \mathrm{d}s$ and $\int_\Gamma(\kappa_j-\kappa_{\max})(u_2(t^n)-u_1(t^n))\phi_{2,j}^{n+1} \mathrm{d}s$ to maintain first order accuracy in time, whose detail is shown as follows
\begin{eqnarray}
&&\int_\Gamma\kappa_j(u_{1,j}(t^{n+1})-u_{2,j}(t^{n+1})){\phi}_{1,j}^{n+1}\mathrm{d}s+\int_\Gamma\kappa_j(u_{2,j}(t^{n+1})-u_{1,j}(t^{n+1})){\phi}_{2,j}^{n+1}\mathrm{d}s\nonumber\\
&&\hspace{5mm}-\int_\Gamma {\kappa_{\max}}(u_{1,j}^{n+1}-u_{2,j}^n) {\phi}_{1,j}^{n+1}\mathrm{d}s
-\int_\Gamma {\kappa_{\max}}(u_{2,j}^{n+1}-u_{1,j}^n) {\phi}_{2,j}^{n+1}\mathrm{d}s\nonumber\\
&&\hspace{5mm}-\int_\Gamma (\kappa_j-{\kappa_{\max}})(u_{1,j}^n-u_{2,j}^n) {\phi}_{1,j}^{n+1} \mathrm{d}s
-\int_\Gamma (\kappa_j-{\kappa_{\max}})(u_{2,j}^n-u_{1,j}^n) {\phi}_{2,j}^{n+1} \mathrm{d}s\nonumber\\
&&=\int_\Gamma\kappa_{\max}({\bf u}_j(t^{n+1})-{\bf u}_j^{n+1}) \cdot {\bf \phi}_j^{n+1}\mathrm{d}s+\int_\Gamma(\kappa_j-\kappa_{\max})({\bf u}_j(t^n)-{\bf u}_j^n) \cdot {\bf\phi}_j^{n+1}\mathrm{d}s\nonumber\\
&&\hspace{5mm}-\int_\Gamma\kappa_j(u_{2,j}(t^n)-u_{2,j}^n)\phi_{1,j}^{n+1}\mathrm{d}s-\int_\Gamma\kappa_j(u_{1,j}(t^n)-u_{1,j}^n)\phi_{2,j}^{n+1}\mathrm{d}s\nonumber\\
&&\hspace{5mm}+\int_\Gamma(\kappa_j-\kappa_{\max})({\bf u}_j(t^{n+1})-{\bf u}_j(t^{n})) \cdot {\bf \phi}_j^{n+1}\mathrm{d}s\nonumber\\
&&\hspace{5mm}-\int_\Gamma\kappa_j({ u}_{2,j}(t^{n+1})-{u}_{2,j}(t^{n})){\phi}_{1,j}^{n+1}\mathrm{d}s-\int_\Gamma\kappa_j({ u}_{1,j}(t^{n+1})-{u}_{1,j}(t^{n})){\phi}_{2,j}^{n+1}\mathrm{d}s.\label{convergence5}
\end{eqnarray}
Substitution of (\ref{convergence5}) for (\ref{convergence4}) and rearranging terms by the error function (\ref{errorfunc}), we can get
\begin{eqnarray}
&&\frac{1}{2\Delta t}(\|{\bf \phi}_j^{n+1}\|^2-\|{\bf \phi}_j^n\|^2)+{\nu}_{\min}\|\nabla {\bf \phi}_j^{n+1}\|^2
+\kappa_{\max}\|{\bf \phi}_{j}^{n+1}\|_\Gamma^2\nonumber\\
&&\le-\frac{1}{\Delta t}({\bf \eta}_j^{n+1}-{\bf \eta}_j^n,{\bf \phi}_j^{n+1})
-({\bf r}_j^{n+1},{\bf \phi}_j^{n+1})-({\nu} \nabla{\bf \eta}_j^{n+1},\nabla{\bf \phi}_j^{n+1})\nonumber\\
&&\hspace{2mm}+\int_\Gamma\kappa_j({ u}_{2,j}(t^{n+1})-{u}_{2,j}(t^{n})){\phi}_{1,j}^{n+1}\mathrm{d}s
+\int_\Gamma\kappa_j({ u}_{1,j}(t^{n+1})-{u}_{1,j}(t^{n})){\phi}_{2,j}^{n+1}\mathrm{d}s\nonumber\\
&&\hspace{2mm}+\int_\Gamma(\kappa_{\max}-\kappa_j)({\bf u}_j(t^{n+1})-{\bf u}_j(t^{n}))\cdot {\bf \phi}_j^{n+1}\mathrm{d}s
              +\int_\Gamma(\kappa_{\max}-\kappa_j)({\bf \eta}_j^n+{\bf \phi}_j^n) \cdot {\bf\phi}_j^{n+1}\mathrm{d}s\nonumber\\
&&\hspace{2mm}+\int_\Gamma\kappa_j(\eta_{2,j}^n+\phi_{2,j}^n)\phi_{1,j}^{n+1}\mathrm{d}s
               +\int_\Gamma\kappa_j(\eta_{1,j}^n+\phi_{1,j}^n)\phi_{2,j}^{n+1}\mathrm{d}s-\int_\Gamma\kappa_{\max} {\bf \eta}_j^{n+1} \cdot {\bf \phi}_j^{n+1}\mathrm{d}s.\label{convergence6}
\end{eqnarray}
Then,  we can bound each term on the RHS of (\ref{convergence6}) with a series of positive constants $\epsilon_m, \ m=1,2,\cdots$. Applying Cauchy-Schwarz, Young's, and Poincar$\acute{e}$ inequalities for the first three terms on the RHS of (\ref{convergence6}), we can yield
\begin{eqnarray}
\frac{1}{\Delta t}({\bf \eta}_j^{n+1}-{\bf \eta}_j^n,{\bf \phi}_j^{n+1})
&\le&\frac{C^2_p}{2\epsilon_1}\left\|\frac{{\bf \eta}_j^{n+1}-{\bf \eta}_j^n}{\Delta t}\right\|^2+\frac{\epsilon_1}{2}\|\nabla {\bf \phi}_j^{n+1} \|^2.\label{convergence7}\\
({\bf r}_j^{n+1},{\bf \phi}_j^{n+1})&\le&\frac{C_p^2}{2\epsilon_2}\|{\bf r}_j^{n+1}\|^2+\frac{\epsilon_2}{2}\|\nabla {\bf \phi}_j^{n+1}\|^2.\label{convergence8}\\
(\nu\nabla{\bf \eta}_j^{n+1},\nabla{\bf \phi}_j^{n+1})
&\le& \frac{ {\nu}_{\max}^2}{2\epsilon_3}\|\nabla {\bf \eta}_j^{n+1}\|^2+\frac{\epsilon_3}{2}\|\nabla {\bf \phi}_j^{n+1}\|^2. \label{convergence9}
\end{eqnarray}
The remaining  seven interface terms on the RHS of (\ref{convergence6}) will be treated in the same way. Utilizing the trace theorem and Young's inequality, we get
\begin{eqnarray}
&&\int_\Gamma\kappa_j({ u}_{2,j}(t^{n+1})-{u}_{2,j}(t^{n})){\phi}_{1,j}^{n+1}\mathrm{d}s+\int_\Gamma\kappa_j({ u}_{1,j}(t^{n+1})-{u}_{1,j}(t^{n})){\phi}_{2,j}^{n+1}\mathrm{d}s\nonumber\\
&&\hspace{3mm}\le  \frac{C_t^4\kappa_j^2}{2\epsilon_4}\|\nabla({\bf u}_j(t^{n+1})-{\bf u}_j(t^{n}))\|^2+\frac{\epsilon_4}{2}\|\nabla{\bf\phi}_j^{n+1}\|^2,\label{convergence10}\\
&&\int_\Gamma(\kappa_{\max}-\kappa_j)({\bf u}_j(t^{n+1})-{\bf u}_j(t^{n}))\cdot{\bf \phi}_j^{n+1}\mathrm{d}s\nonumber\\
&&\hspace{3mm}\le \frac{C_t^4(\kappa_{\max}-\kappa_j)^2}{2\epsilon_5}\|\nabla({\bf u}_j(t^{n+1})-{\bf u}_j(t^{n}))\|^2+\frac{\epsilon_5}{2}\|\nabla{\bf \phi}_j^{n+1}\|^2,\label{convergence11}\\
&&\int_\Gamma(\kappa_{\max}-\kappa_j)({\bf \eta}_j^n+{\bf \phi}_j^n)\cdot{\bf\phi}_j^{n+1}\mathrm{d}s\nonumber\\
&&\hspace{3mm}\le\frac{\kappa_{\max}-\kappa_j}{2}\|{\bf \phi}_j^n\|_\Gamma^2+\frac{\kappa_{\max}-\kappa_j}{2}\|{\bf \phi}_j^{n+1}\|_\Gamma^2+\frac{C_t^4(\kappa_{\max}-\kappa_j)^2}{2\epsilon_6}\|\nabla {\bf \eta}_j^n\|^2+\frac{\epsilon_6}{2}\|\nabla{\bf \phi}_j^{n+1}\|^2,\hspace{10mm}\label{convergence12}\\
&&\int_\Gamma\kappa_j(\eta_{2,j}^n+\phi_{2,j}^n)\phi_{1,j}^{n+1}\mathrm{d}s+\int_\Gamma\kappa_j(\eta_{1,j}^n+\phi_{1,j}^n)\phi_{2,j}^{n+1}\mathrm{d}s\nonumber\\
&&\hspace{3mm}\le\frac{\kappa_j}{2}\|{\bf \phi}_j^n\|_\Gamma^2+\frac{\kappa_j}{2}\|{\bf \phi}_j^{n+1}\|_\Gamma^2+\frac{C_t^4\kappa_j^2}{2\epsilon_7}\|\nabla{\bf \eta}_j^n\|^2+\frac{\epsilon_7}{2}\|\nabla{\bf \phi}_j^{n+1}\|^2,\label{convergence13}\\
&&\int_\Gamma\kappa_{\max} {\bf \eta}_j^{n+1}\cdot{\bf \phi}_j^{n+1}\mathrm{d}s\le \frac{C_t^4\kappa_{\max}^2}{2\epsilon_8}\|\nabla{\bf \eta}_j^{n+1}\|^2+\frac{\epsilon_8}{2}\|\nabla{\bf \phi}_j^{n+1}\|^2.\label{convergence14}
\end{eqnarray}
Setting $\epsilon_1=\cdots=\epsilon_8=\frac{\nu_{\min}}{8}$ and substituting the above inequalities (\ref{convergence7})-(\ref{convergence14}) into (\ref{convergence6}), we can conclude
\begin{eqnarray}
&&\frac{1}{2\Delta t}(\|{\bf \phi}_j^{n+1}\|^2-\|{\bf \phi}_j^{n}\|^2)+\frac{\nu_{\min}}{2}\|\nabla {\bf \phi}_j^{n+1}\|^2
+\frac{\kappa_{\max}}{2}(\|{\bf\phi}_j^{n+1}\|_\Gamma^2-\|{\bf \phi}_j^{n}\|_\Gamma^2)\nonumber\\
&&\hspace{5mm}\le C\bigg\{\max(\kappa_{\max}^2,\nu^2_{\max})\|\nabla{\bf \eta}_j^{n+1}\|^2+\max(\kappa_j^2,(\kappa_{\max}-\kappa_j)^2)\|\nabla{\bf \eta}_j^{n}\|^2\nonumber\\
&&\hspace{10mm}+\max(\kappa_j^2,(\kappa_{\max}-\kappa_j)^2)\|\nabla({\bf u}_j(t^{n+1})-{\bf u}_j(t^{n}))\|^2\nonumber\\
&&\hspace{10mm}+ \left\|\frac{{\bf \eta}_j^{n+1}-{\bf \eta}_j^n}{\Delta t}\right\|^2+\|{\bf r}_j^{n+1}\|^2\bigg\}.\label{convergence15}
\end{eqnarray}
Multiplying through (\ref{convergence15}) by $2\Delta t$ and summing over $n=0,1,\cdots,N-1 $ it follows
\begin{eqnarray}
&&\|{\bf \phi}_j^{N}\|^2+\kappa_{\max}\Delta t\|{\bf \phi}_j^{N}\|_\Gamma^2+{\nu}_{\min}\Delta t\sum\limits_{n=0}^{N-1}\|\nabla{\bf \phi}_j^{n+1}\|^2\nonumber\\
&&\hspace{5mm}\le \|{\bf \phi}_j^0\|^2+\kappa_{\max}\Delta t\|{\bf \phi}_j^0\|_\Gamma^2+C\Delta t\sum\limits_{n=0}^{N-1}
\bigg\{\|\nabla{\bf \eta}_j^{n+1}\|^2+\|\nabla{\bf \eta}_j^{n}\|^2\nonumber\\
&&\hspace{8mm}+\|\nabla({\bf u}_j(t^{n+1})-{\bf u}_j(t^{n}))\|^2+\left\|\frac{{\bf \eta}_j^{n+1}-{\bf \eta}_j^n}{\Delta t}\right\|^2
+\|{\bf r}_j^{n+1}\|^2\bigg\}.\label{convergence16}
\end{eqnarray}
For the last three terms in (\ref{convergence16}), we have the following estimates
\begin{eqnarray}
&&\Delta t\sum\limits_{n=0}^{N-1} \left\|\frac{{\bf\eta}_j^{n+1}-{\bf \eta}_j^n}{\Delta t}\right\|^2\le \int_0^{t^{n+1}}\|{\bf \eta}_{j,t}\|^2\mathrm{d}t
\le\|{\bf \eta}_{j,t}\|_{L^2(0,T;L^2(\Omega))}^2, \label{convergence171} \\
&&\Delta t\sum\limits_{n=0}^{N-1} \|\nabla({\bf u}_j(t^{n+1})-{\bf u}_j(t^{n}))\|^2 \le {\Delta t}^2\int_0^{t^{n+1}}\|\nabla{\bf u}_{j,t}\|^2\mathrm{d}t
\le{\Delta t}^2\|\nabla{\bf u}_{j,t}\|_{L^2(0,T;L^2(\Omega))}^2, \hspace{10mm}\label{convergence172}\\
&&\Delta t\sum\limits_{n=0}^{N-1}\|{\bf r}_j^{n+1}\|^2\le {\Delta t}^2\int_0^{t^{n+1}}\|{\bf u}_{j,tt}\|^2\mathrm{d}t\le{\Delta t}^2\|{\bf u}_{j,tt}\|_{L^2(0,T;L^2(\Omega))}^2.\label{convergence173}
\end{eqnarray}
Taking infimum over ${\bf v}_j^n\in X_h$, using the triangle inequality to $\|{\bf \phi}_j^0\|^2+\kappa_{\max}\Delta t\|{\bf \phi}_j^0\|_\Gamma^2$, and combining with the inequalities  (\ref{convergence171})-(\ref{convergence173}), we can yield
\begin{eqnarray}\label{convergencef}
&&\|{\bf \phi}_j^{N}\|^2+\kappa_{\max}\Delta t\|{\bf \phi}_j^{N}\|_\Gamma^2+{\nu}_{\min}\Delta t\sum\limits_{n=0}^{N-1}\|\nabla{\bf \phi}_j^{n+1}\|^2\nonumber\\
&&\le C\bigg\{\|{\bf u}_j(0)-{\bf u}_j^0\|^2+\kappa_{\max}\Delta t\|{\bf u}_j(0)-{\bf u}_j^0\|_\Gamma^2+{\Delta t}^2\|\nabla{\bf u}_{j,t}\|_{L^2(0,T;L^2{\Omega})}^2\bigg.\nonumber\\
&&\hspace{8mm}\bigg.+{\Delta t}^2\|{\bf u}_{j,tt}\|_{L^2(0,T;L^2(\Omega))}^2+\inf\limits_{{\bf v}_j^0\in X_h}\left\{\|{\bf \eta}_j^0\|^2+\kappa_{\max}\Delta t\|{\bf \eta}_j^0\|_\Gamma^2\right\}\bigg.\nonumber\\
&&\hspace{8mm}\bigg. + \inf\limits_{{\bf v}_j\in X_h}\|{\bf \eta}_{j,t}\|_{L^2(0,T;L^2(\Omega))}^2 + \Delta t \inf\limits_{{\bf v}_j^n\in X_h} \sum\limits_{n=0}^{N-1} \left\{\|\nabla{\bf \eta}_j^{n+1}\|^2+\|\nabla{\bf \eta}_j^{n}\|^2\right\}\bigg\},
\end{eqnarray}
where $C$ is a generic constant independent of the mesh size $h$ and time step $\Delta t$. Then, we can further deal with the last them of (\ref{convergencef}) as follows
\begin{eqnarray*}
\Delta t \inf\limits_{{\bf v}_j^n\in X_h} \sum\limits_{n=0}^{N-1} \left\{\|\nabla{\bf \eta}_j^{n+1}\|^2+\|\nabla{\bf \eta}_j^{n}\|^2\right\}\le 2T \max\limits_{n=0,1,\cdots,N} \inf\limits_{{\bf v}_j^n\in X_h}\|\nabla{\bf\eta}_j^n\|^2.
\end{eqnarray*}
By applying triangle inequality and rearranging constants, we can summarize the final result (\ref{A1convergence}).
\end{proof}

\begin{rem} \label{remk}
Note that the existing ensemble algorithms in \cite{wang zhu, random PDEs2, ensemble1, ensemble2, ensemble3, ensemble4, ensemble5} both need to assume a small perturbation constraint to the random parameters. However, because of the particularity of the parameter that the friction parameter $\kappa$ is a random function independent of time and space, the A1 algorithm can well avoid this constraint in the above analysis, which is also mentioned in Section 4 of \cite{ensemble2}.
Moreover, to further indicate this advantage, we test the $\kappa$ from $10^{-2}$ to $10^1$ in the numerical tests, which is usually between $10^{-3}$ and $10^3$ \cite{heat-linear}.
\end{rem}

\section{The ensemble scheme for CASE 2}
Based on the analytical result of CASE 1, we extend the method to the heat-heat coupled model with three random coefficients $\kappa, \nu_1$, and $\nu_2$. Then, we can propose the following ensemble algorithm for CASE2.

\textbf{{{Ensemble Algorithm 2 (A2)}}}

Let $\Delta t>0$, $\bar\nu_i^n=\frac{1}{J}\sum\limits_{j=1}^J\nu_{i,j}({\bf x},t^n)$, ${\bf u}_j^0\in X_h$ and ${\bf f}_j\in L^2(H^{-1}(\Omega); 0,T))$ for $j=1,\cdots,J$. Given ${\bf u}_j^n\in X_h$, find ${\bf u}_j^{n+1}\in X_h$ for $n=0,1,2,\ldots,N-1$ and $j=1,\cdots,J$,
\begin{eqnarray}
&&(\frac{u_{1,j}^{n+1}-u_{1,j}^{n}}{\Delta t},v_1)+({\bar\nu}_1^{n+1}\nabla u_{1,j}^{n+1},\nabla v_1)+((\nu_{1,j}^{n+1}-{\bar\nu}_1^{n+1})\nabla u_{1,j}^n,\nabla v_1)\nonumber\\
&&\hspace{10mm}+\int_\Gamma {\kappa_{\max}}(u_{1,j}^{n+1}-u_{2,j}^n) v_1\mathrm{d}s+\int_\Gamma (\kappa_j-{\kappa_{max}})(u_{1,j}^n-u_{2,j}^n) v_1 \mathrm{d}s\nonumber\\
&&\hspace{65mm}=(f_{1,j}^{n+1},v_1) \ \ \ \ \forall v_1\in X_{1,h},\label{algorithm2.1}\\
&&(\frac{u_{2,j}^{n+1}-u_{2,j}^{n}}{\Delta t},v_2)+({\bar\nu}_2^{n+1}\nabla u_{2,j}^{n+1},\nabla v_2)+((\nu_{2,j}^{n+1}-{\bar\nu}_2^{n+1})\nabla u_{2,j}^n,\nabla v_2)\nonumber\\
&&\hspace{10mm}+\int_\Gamma {\kappa_{\max}}(u_{2,j}^{n+1}-u_{1,j}^n) v_2\mathrm{d}s+\int_\Gamma (\kappa_j-{\kappa_{\max}})(u_{2,j}^n-u_{1,j}^n) v_2 \mathrm{d}s\nonumber\\
&&\hspace{65mm}=(f_{2,j}^{n+1},v_2) \ \ \ \  \forall v_2\in X_{2,h}.\label{algorithm2.2}
\end{eqnarray}

Suppose the following two conditions are valid:

(\romannumeral3) There exists a positive constant $\theta$ such that, for any $t \in [0,T]$,
\begin{eqnarray}
\min\limits_{{\bf x}\in \bar{\Omega}} \bar\nu_i({\bf x},t)\ge \theta.\label{nuj}
\end{eqnarray}

(\romannumeral4)There exist positive constants $\theta_-$ and $\theta_+$ such that, for any $t\in [0,T]$,
\begin{eqnarray}
\theta_-\le |\nu_{i,j}({\bf x},t)-\bar{\nu}_i({\bf x},t)|_{\infty}\le\theta_+.\label{nujbarnu}
\end{eqnarray}
\begin{theorem}(A2 Stability)
Suppose that ${\bf f}_j\in L^2(H^{-1}(\Omega);0,T)$, ${\bf u}_j^{n+1}\in X_h$ satisfy (\ref{algorithm2.1})-(\ref{algorithm2.2}) for each $n\in \{0,1,2,\cdots, N-1\}$, and the conditions (\romannumeral3)-(\romannumeral4) hold, the algorithm A2 is stable on the premise that
\begin{eqnarray}
\theta>\theta_+. \label{stacondition}
\end{eqnarray}
Then there exists a generic positive constant C independent of $h, \Delta t, J$, such that the numerical solution ${\bf u}_j^{n+1}$ satisfies
\begin{eqnarray}
\|{\bf u}_j^{N}\|^2+\theta_-\Delta t\|\nabla {\bf u}_j^{N}\|^2+(\theta-\theta_+)\Delta t\sum\limits_{n=0}^{N-1}\|\nabla{\bf u}_j^{n+1}\|^2
+\kappa_{\max}\Delta t\|{\bf u}_j^{N}\|^2_\Gamma \nonumber\\
\le \|{\bf u}_j^0\|^2+C\Delta t\|\nabla {\bf u}_j^0\|^2+C\Delta t\|{\bf u}_j^0\|_\Gamma^2+C\Delta t\sum\limits_{n=0}^{N-1}\|{\bf f}_j^{n+1}\|_{-1}^2.
\end{eqnarray}
\end{theorem}

\begin{proof}
We can get a detailed proof process by combining Theorem 1 of \cite{wang zhu} and Theorem 3.1 of this paper.
\end{proof}

\begin{theorem}(A2 Error Estimate)
Let ${\bf u}_j(t^{n+1})$ and ${\bf u}_j^{n+1}$ be the solutions of (\ref{Moweakform}) and the algorithm A2 at time $t^{n+1}$, respectively. Assume conditions (\romannumeral3)-(\romannumeral4) and the stability condition (\ref{stacondition}) hold. Then, there exists a generic positive constant C independent of $h, \Delta t, J$ for any $n\in \{0,1,2,\cdots, N-1\}$ such that
\begin{eqnarray}
&&\|{\bf u}_j(t^{N})-{\bf u}_j^{N}\|^2+\kappa_{\max}\Delta t\|{\bf u}_j(t^{N})-{\bf u}_j^{N}\|_\Gamma^2\nonumber\\
&&\hspace{5mm}+\theta_-\Delta t\|\nabla{\bf u}_j(t^{N})-\nabla{\bf u}_j^{N}\|^2+(\theta-\theta_+)\Delta t\sum\limits_{n=0}^{N-1}\|\nabla{\bf u}_j(t^{n+1})-\nabla{\bf u}_j^{n+1}\|^2\nonumber\\
&&\le C\bigg\{\|{\bf u}_j(0)-{\bf u}_j^0\|^2+\kappa_{\max}\Delta t\|{\bf u}_j(0)-{\bf u}_j^0\|_\Gamma^2+{\Delta t}^2\|\nabla{\bf u}_{j,t}\|_{L^2(0,T;L^2{\Omega})}^2+{\Delta t}^2\|{\bf u}_{j,tt}\|_{L^2(0,T;L^2(\Omega))}^2\bigg.\nonumber\\
&&\hspace{5mm}\bigg.+\inf\limits_{{\bf v}_j^0\in X_h}\left\{\|{\bf u}_j(0)-{\bf v}_j^0\|^2+\kappa_{\max}\Delta t\|\nabla({\bf u}_j(0)-{\bf v}_j^0)\|_\Gamma^2\right\}+\inf\limits_{{\bf v}_j\in X_h}\|({\bf u}_j(0)-{\bf v}_j)_t\|_{L^2(0,T;L^2(\Omega))}^2 \bigg.\nonumber\\
&&\hspace{5mm}\bigg.+T \max\limits_{n=0,1,\cdots,N}\inf\limits_{{\bf v}_j^n\in X_h}\|\nabla({\bf u}_j(t^n)-{\bf v}_j^n)\|^2 \bigg\}.\label{error2.0}
\end{eqnarray}
\end{theorem}
\begin{proof}
Restricting test function ${\bf v}$ to $X_h$ and subtracting (\ref{algorithm2.1})-(\ref{algorithm2.2}) from (\ref{Moweakform}), we can derive the following error equation by the similar way of (\ref{convergence2}).
\begin{eqnarray}
&&({\bf r}^{n+1},{\bf v})+\left(\frac{{\bf u}_j(t^{n+1})-{\bf u}_j^{n+1}}{\Delta t}-\frac{{\bf u}_j(t^{n})-{\bf u}_j^{n}}{\Delta t},{\bf v}\right)+({\bar\nu}^{n+1}(\nabla{\bf u}_j(t^{n+1})-\nabla{\bf u}_j^{n+1}),\nabla{\bf v}) \nonumber\\
&&+((\nu_j^{n+1}-{\bar\nu}^{n+1})(\nabla{\bf u}_j(t^{n})-\nabla{\bf u}_j^{n}),\nabla{\bf v})
+((\nu_j^{n+1}-{\bar\nu}^{n+1})(\nabla{\bf u}_j(t^{n+1})-\nabla{\bf u}_j(t^n)),\nabla{\bf v})\nonumber\\
&&+\int_\Gamma\kappa_j(u_{1,j}(t^{n+1})-u_{2,j}(t^{n+1})){v_1}\mathrm{d}s +\int_\Gamma\kappa_j(u_{2,j}(t^{n+1})-u_{1,j}(t^{n+1})){v_2}\mathrm{d}s\nonumber\\
&&-\int_\Gamma {\kappa_{\max}}(u_{1,j}^{n+1}-u_{2,j}^n) v_1\mathrm{d}s
-\int_\Gamma {\kappa_{\max}}(u_{2,j}^{n+1}-u_{1,j}^n) v_2\mathrm{d}s\nonumber\\
&&-\int_\Gamma (\kappa_j-{\kappa_{\max}})(u_{1,j}^n-u_{2,j}^n) v_1 \mathrm{d}s
-\int_\Gamma (\kappa_j-{\kappa_{\max}})(u_{2,j}^n-u_{1,j}^n) v_2\mathrm{d}s=0.\label{error2.1}
\end{eqnarray}
Choosing ${\bf v}={\bf \phi}_j^{n+1}$, we obtain
\begin{eqnarray}
&&\frac{1}{2\Delta t}(\|{\bf \phi}_j^{n+1}\|^2-\|{\bf \phi}_j^n\|^2)+\theta\|\nabla {\bf \phi}_j^{n+1}\|^2
+\kappa_{\max}\|{\bf \phi}_{j}^{n+1}\|_\Gamma^2\nonumber\\
&&\le-\frac{1}{\Delta t}({\bf \eta}_j^{n+1}-{\bf \eta}_j^n,{\bf \phi}_j^{n+1})-({\bf r}_j^{n+1},{\bf \phi}_j^{n+1})-\mathop{\uwave{((\nu_j^{n+1}-\bar{\nu}^{n+1})\nabla{\bf\phi}_j^{n},\nabla{\bf\phi}_j^{n+1})}}\limits_{\textcircled{1}}
-\mathop{\uwave{(\bar{\nu}^{n+1}\nabla{\bf\eta}_j^{n+1},\nabla{\bf\phi}_j^{n+1})}}\limits_{\textcircled{2}}\nonumber\\
&&\hspace{4mm}-\mathop{\uwave{(({\nu}_j^{n+1}-\bar{\nu}^{n+1})\nabla{\bf\eta}^n_j,\nabla{\bf\phi}_j^{n+1})}}\limits_{\textcircled{3}}
-\mathop{\uwave{((\nu_j^{n+1}-\bar{\nu}^{n+1})(\nabla{\bf u}_j(t^{n+1})-\nabla{\bf u}_j(t^n)),\nabla{\bf\phi}_j^{n+1})}}\limits_{\textcircled{4}}\nonumber\\
&&\hspace{4mm}+\int_\Gamma\kappa_j({ u}_{2,j}(t^{n+1})-{u}_{2,j}(t^{n})){\phi}_{1,j}^{n+1}\mathrm{d}s
+\int_\Gamma\kappa_j({ u}_{1,j}(t^{n+1})-{u}_{1,j}(t^{n})){\phi}_{2,j}^{n+1}\mathrm{d}s\nonumber\\
&&\hspace{4mm}+\int_\Gamma(\kappa_{\max}-\kappa_j)({\bf u}_j(t^{n+1})-{\bf u}_j(t^{n}))\cdot{\bf \phi}_j^{n+1}\mathrm{d}s+\int_\Gamma(\kappa_{\max}-\kappa_j)({\bf u}_j(t^n)-{\bf u}_j^n)\cdot{\bf\phi}_j^{n+1}\mathrm{d}s\nonumber\\
&&\hspace{4mm}+\int_\Gamma\kappa_j(u_{2,j}(t^n)-u_{2,j}^n)\phi_{1,j}^{n+1}\mathrm{d}s
+\int_\Gamma\kappa_j(u_{1,j}(t^n)-u_{1,j}^n)\phi_{2,j}^{n+1}\mathrm{d}s-\int_\Gamma\kappa_{\max} {\bf \eta}_j^{n+1}\cdot{\bf \phi}_j^{n+1}\mathrm{d}s.
\label{error2.2}
\end{eqnarray}
Next, we need to bound the terms on the RHS of (\ref{error2.2}). The main difference between the proof of Theorem 3.2 lies in the estimates of \textcircled{1}-\textcircled{4}, which requires more rigorous analysis to replace the inequality (\ref{convergence9}). In addition, this also means that the parameter $\epsilon_3$ does not exist in this demonstration. For the estimation of \textcircled{1}, \textcircled{3}, \textcircled{4}, we can directly refer to the analysis in Theorem 3 of \cite{wang zhu} (Page 865) with four positive constants $\beta_s (s=0,\cdots, 3)$.
Meanwhile, for \textcircled{2}, we have the following inequality
\begin{eqnarray*}
(\bar{\nu}^{n+1}\nabla{\bf\eta}_j^{n+1},\nabla{\bf\phi}_j^{n+1})\le|\bar{\nu}^{n+1}|_\infty\left(\frac{\|\nabla{\bf\eta}_j^{n+1}\|^2}{2\epsilon_9}
+\frac{\epsilon_9\|\nabla{\bf\phi}_j^{n+1}\|^2}{2}\right).
\end{eqnarray*}
Compared with \cite{wang zhu}, in order to ensure that no conditions are imposed or strengthened, we make the same estimation for $({\bf{r}}_j^{n+1}, {\bf\phi}_j^{n+1})$ as \cite{wang zhu} (Page 865) to replace inequality (\ref{convergence8}), which also means that the parameter $\epsilon_3$ does not appear. Therefore, we can choose $\epsilon_1=\epsilon_4=\epsilon_5=\epsilon_6=\epsilon_7=\epsilon_8=\frac{|\bar{\nu}^{n+1}|_\infty\beta_0}{7}$ in Theorem 3.2 and $\epsilon_9=\frac{\beta_0}{7}$ and combine 
the analysis of Theorem 3.2 and Theorem 3 of \cite{wang zhu} to arrive
\begin{eqnarray}
&&\frac{1}{2\Delta t}(\|{\bf \phi}_j^{n+1}\|^2-\|{\bf \phi}_j^{n}\|^2)+\frac{|\nu_j^{n+1}-\bar\nu^{n+1}|_\infty}{2}(\|\nabla {\bf \phi}_j^{n+1}\|^2-\|\nabla {\bf \phi}_j^{n}\|^2)+\frac{\kappa_{\max}}{2}(\|{\bf \phi}_j^{n+1}\|_\Gamma^2-\|{\bf \phi}_j^{n}\|_\Gamma^2)\nonumber\\
&&\hspace{10mm}+\bigg(\theta-\frac{\beta_0}{2}|\bar\nu^{n+1}|_\infty-\frac{\beta_1+\beta_2+2}{2}|\nu_j^{n+1}-\bar\nu^{n+1}|_\infty-\frac{\beta_3}{2}\bigg)
\|\nabla{\bf \phi}_j^{n+1}\|^2\nonumber\\
&&\le\frac{|\bar\nu^{n+1}|_\infty}{2\beta_0}\|\nabla{\bf\eta}_j^{n+1}\|^2
+\frac{|\nu_j^{n+1}-\bar\nu^{n+1}|_\infty}{2\beta_1}\|\nabla{\bf\eta}_j^n\|^2\nonumber\\
&&\hspace{4mm}+\frac{|\nu_j^{n+1}-\bar\nu^{n+1}|_\infty}{2\beta_2}\|\nabla({\bf u}_j(t^{n+1})-{\bf u}_j(t^{n}))\|^2 +\frac{\|{\bf r}_j^{n+1}\|_{-1}^2}{2\beta_3}\nonumber\\
&&\hspace{4mm}
+C\bigg\{\frac{\kappa_{\max}^2}{{|\bar\nu^{n+1}|_\infty\beta_0}}\|\nabla{\bf \eta}_j^{n+1}\|^2
+\max\left(\frac{\kappa_j^2}{|\bar\nu^{n+1}|_\infty\beta_0},\frac{(\kappa_{\max}-\kappa_j)^2}{|\bar\nu^{n+1}|_\infty\beta_0}\right)\|\nabla({\bf u}_j(t^{n+1})-{\bf u}_j(t^{n}))\|^2\nonumber\\
&&\hspace{12mm} +\max\left(\frac{\kappa_j^2}{{|\bar\nu^{n+1}|_\infty\beta_0}},\frac{(\kappa_{\max}-\kappa_j)^2}{{|\bar\nu^{n+1}|_\infty\beta_0}}\right)\|\nabla{\bf \eta}_j^{n}\|^2
+ \frac{1}{{|\bar\nu^{n+1}|_\infty\beta_0}}\left\|\frac{{\bf \eta}_j^{n+1}-{\bf \eta}_j^n}{\Delta t}\right\|^2\bigg\}.\hspace{10mm}\label{error2.3}
\end{eqnarray}
We can substitute the selection $\beta_0=\frac{\delta|\nu_j^{n+1}-\bar\nu^{n+1}|_\infty}{2|\bar\nu_{n+1}|_\infty}$, $\beta_1=\beta_2=\frac{\delta}{2}$, and $\beta_3=\frac{\delta|\nu_j^{n+1}-\bar\nu^{n+1}|_\infty}{2}$ with a positive constant $\delta$ in Reference \cite{wang zhu} into (\ref{error2.3}) to yield
\begin{eqnarray}
&&\frac{1}{2\Delta t}(\|{\bf \phi}_j^{n+1}\|^2-\|{\bf \phi}_j^{n}\|^2)+\frac{|\nu_j^{n+1}-\bar\nu^{n+1}|_\infty}{2}(\|\nabla {\bf \phi}_j^{n+1}\|^2-\|\nabla {\bf \phi}_j^{n}\|^2)\nonumber\\
&&+\bigg(\theta-(1+\delta)|\nu_j^{n+1}-\bar\nu^{n+1}|_\infty\bigg)
\|\nabla{\bf \phi}_j^{n+1}\|^2
+\frac{\kappa_{\max}}{2}(\|{\bf \phi}_j^{n+1}\|_\Gamma^2-\|{\bf \phi}_j^{n}\|_\Gamma^2)\nonumber\\
&&\hspace{1mm}\le C\bigg\{\max\bigg(\frac{|\bar\nu^{n+1}|_\infty}{\delta|\nu_j^{n+1}-\bar\nu^{n+1}|_\infty},\frac{\kappa_{\max}^2}{\delta|\nu_j^{n+1}-\bar\nu^{n+1}|_\infty}\bigg)
\|\nabla{\bf\eta}_j^{n+1}\|^2\nonumber\\
&&\hspace{4mm}+\max\bigg(\frac{|\nu_j^{n+1}-\bar\nu^{n+1}|_\infty}{\delta},\frac{\kappa_j^2}{\delta|\nu_j^{n+1}-\bar\nu^{n+1}|_\infty},
\frac{(\kappa_{\max}-\kappa_j)^2}{\delta|\nu_j^{n+1}-\bar\nu^{n+1}|_\infty}\bigg)\|\nabla{\bf\eta}_j^n\|^2\nonumber\\
&&\hspace{4mm}+\max\bigg(\frac{|\nu_j^{n+1}-\bar\nu^{n+1}|_\infty}{\delta},\frac{\kappa_j^2}{\delta|\nu_j^{n+1}-\bar\nu^{n+1}|_\infty},
\frac{(\kappa_{\max}-\kappa_j)^2}{\delta|\nu_j^{n+1}-\bar\nu^{n+1}|_\infty}\bigg)\|\nabla({\bf u}_j(t^{n+1})-{\bf u}_j(t^{n}))\|^2\nonumber\\
&&\hspace{4mm}+\frac{\|{\bf r}_j^{n+1}\|_{-1}^2}{\delta|\nu_j^{n+1}-\bar\nu^{n+1}|_\infty}
+\frac{1}{\delta|\nu_j^{n+1}-\bar\nu^{n+1}|_\infty}\left\|\frac{{\bf \eta}_j^{n+1}-{\bf \eta}_j^n}{\Delta t}\right\|^2\bigg\}.\label{error2.4}
\end{eqnarray}
Setting $\delta=\frac{\theta-\theta_+}{2\theta_+}$, we have $\theta-(1+\delta)|\nu_j^{n+1}-\bar\nu^{n+1}|_\infty>\frac{\theta-\theta_+}{2}>0$ based on the stability condition (\ref{stacondition}) and the upper bound in condition (\ref{nujbarnu}). The rest of the detailed analysis process is similar to Theorem 3.2, and the final result can be derived as follows
\begin{eqnarray}
&&\|{\bf \phi}_j^{N}\|^2+\kappa_{\max}\Delta t\|{\bf \phi}_j^{N}\|_\Gamma^2+(\theta-\theta_+)\Delta t\sum\limits_{n=0}^{N-1}\|\nabla{\bf\phi}_j^{n+1}\|^2+\theta_-\Delta t\|\nabla{\bf \phi}_j^{N}\|^2\nonumber\\
&&\le C\bigg\{\|{\bf u}_j(0)-{\bf u}_j^0\|^2+\kappa_{\max}\Delta t\|{\bf u}_j(0)-{\bf u}_j^0\|_\Gamma^2+{\Delta t}^2\|\nabla{\bf u}_{j,t}\|_{L^2(0,T;L^2{\Omega})}^2\bigg.\nonumber\\
&&\hspace{8mm}\bigg.+{\Delta t}^2\|{\bf u}_{j,tt}\|_{L^2(0,T;L^2(\Omega))}^2+\inf\limits_{{\bf v}_j^0\in X_h}\left\{\|{\bf \eta}_j^0\|^2+\kappa_{\max}\Delta t\|{\bf \eta}_j^0\|_\Gamma^2\right\}+ \inf\limits_{{\bf v}_j\in X_h}\|{\bf \eta}_{j,t}\|_{L^2(0,T;L^2(\Omega))}^2\bigg.\nonumber\\
&&\hspace{8mm}\bigg.+\Delta t \inf\limits_{{\bf v}_j^n\in X_h} \sum\limits_{n=0}^{N-1} \left\{\|\nabla{\bf \eta}_j^{n+1}\|^2+\|\nabla{\bf \eta}_j^{n}\|^2\right\}\bigg\},\label{error2.5}
\end{eqnarray}
where $C$ is a generic positive constant independent of $h, \Delta t$, and $J$. By the triangle inequality, we have the final result (\ref{error2.0}).
\end{proof}
\begin{rem}
We proposed the A2 for the heat-heat coupled problem with three random coefficients, which are more complex than the single domain random problem \cite{wang zhu}. Fortunately, compared with the results of \cite{wang zhu}, we do not need to strengthen or impose any constraint conditions to obtain unconditional stability and convergence.
\end{rem}

Note that we calculate the expectation of diffusion coefficient $\bar\nu_i^n$ at each time step as the key point of the ensemble idea in the A2 algorithm. Interestingly, we can further optimize the A2 algorithm. On the basis of calculating the expectation of $\nu_{i,j}(\bf{x},t^n)$, we also calculate the average in time. In this way, several linear equations we get not only share a coefficient matrix in the sample, but also in time, so as to further reduce storage requirements and save computing costs. The optimized algorithm is shown below

\textbf{{{Ensemble Algorithm 3} (A3)}}

Let $\Delta t>0$, $\bar\nu_i=\frac{1}{N}\sum\limits_{n=1}^N\bar\nu_i^n$, ${\bf u}_j^0\in X_h$ and ${\bf f}_j\in L^2(H^{-1}(\Omega); 0,T))$ for $j=1,\cdots,J$. Given ${\bf u}_j^n\in X_h$, find ${\bf u}_j^{n+1}\in X_h$ for $n=0,1,2,\cdots,N-1$ and $j=1,\cdots,J$,
\begin{eqnarray}
&&(\frac{u_{1,j}^{n+1}-u_{1,j}^{n}}{\Delta t},v_1)+(\bar\nu_1 \nabla u_{1,j}^{n+1},\nabla v_1)+((\nu_{1,j}^{n+1}-\bar\nu_1)\nabla u_{1,j}^n,\nabla v_1)\nonumber\\
&&\hspace{3mm}+\int_\Gamma {\kappa_{\max}}(u_{1,j}^{n+1}-u_{2,j}^n) v_1\mathrm{d}s+\int_\Gamma (\kappa_j-{\kappa_{\max}})(u_{1,j}^n-u_{2,j}^n) v_1 \mathrm{d}s\nonumber\\
&&\hspace{65mm}=(f_{1,j}^{n+1},v_1) \ \ \ \ \forall v_1\in X_{1,h},\label{algorithm3.1}\\
&&(\frac{u_{2,j}^{n+1}-u_{2,j}^{n}}{\Delta t},v_2)+(\bar\nu_2 \nabla u_{2,j}^{n+1},\nabla v_2)+(\nu_{2,j}^{n+1}-\bar\nu_2)\nabla u_{2,j}^n,\nabla v_2)\nonumber\\
&&+\int_\Gamma {\kappa_{\max}}(u_{2,j}^{n+1}-u_{1,j}^n) v_2\mathrm{d}s+\int_\Gamma (\kappa_j-{\kappa_{\max}})(u_{2,j}^n-u_{1,j}^n) v_2 \mathrm{d}s\nonumber\\
&&\hspace{65mm}=(f_{2,j}^{n+1},v_2) \ \ \ \  \forall v_2\in X_{2,h}.\label{algorithm3.2}
\end{eqnarray}

The proof for the stability and the convergence of A3 is similar to A2, and the parameter condition $\theta>\theta_+$ is also required for A3.

\begin{rem}
If it is easy to identify the time and space maximum of the coefficient $\nu_{i,j}({\bf x},t)$, the parameter conditions $\theta>\theta_+$ of A2 and A3 is not necessary, and the detailed analysis process is similar to CASE 1. In real, selecting the maximum of the coefficient $\nu_{i,j}({\bf x},t)$ is not easy work and may take a lot of time. Hence, the usual method is to calculate the expectation of the diffusion coefficient $\nu_{i,j}({\bf x},t^n)$ with a small perturbation constraint. This also gives us an idea that if the random parameter has a large perturbation, we can compute the maximum value or $L^{\infty}$ norm of the random parameter as the key of the ensemble method.
\end{rem}

\section{Numerical results}
In this section, we present three numerical tests to illustrate the approximate accuracy and efficiency of the proposed ensemble algorithm for the heat-heat coupled system with three random coefficients. In the first numerical experiment, we test a smooth problem to verify the convergence of the three ensemble algorithm. Meanwhile, we also record the elapsed CPU time to check the superiority of our algorithms. In the second example, we test the energy attenuation results of the three algorithms to prove that the proposed algorithms are unconditionally stable. Finally, we design a steel-titanium composite plate fuel cell 3D model to observe its internal heat conduction phenomenon. In both tests, the finite element spaces are constructed by linear Lagrangian elements (P1) for the heat-heat coupled model. All the numerical experiments are implemented by the open software FreeFEM++ \cite{FEM++}.

\subsection{Smooth problem}
In the first numerical experiment, we test a smooth problem with an exact solution adapted from \cite{heat-linear} to check the convergence rate of our ensemble algorithms. Assume $\Omega_1=[0,1]\times[0,1]$ and $\Omega_2=[0,1]\times[-1,0]$, the interface $\Gamma$ of the current
computational domain is the portion of the $x$-axis from 0 to 1, $\hat{n}_1=[0,-1]^{T}$ and $\hat{n}_2=[0,1]^T$. The exact solution is selected as follows:
\begin{eqnarray*}
&&u_1(t,x,y)=ax(1-x)(1-y)e^{-t},\\
&&u_2(t,x,y)=ax(1-x)(c_1+c_2y+c_3 y^2)e^{-t},
\end{eqnarray*}
which also determines the Dirichlet boundary condition, initial condition, and source terms of this smooth problem. The selection of above constants $c_1, c_2, c_3$ should be determined by (\ref{heat1})-(\ref{heat4}). Inspired by \cite{heat-linear}, we can choose 
\begin{eqnarray*}
c_1=1+\frac{\nu_1}{\kappa}, c_2=\frac{-\nu_1}{\nu_2}, c_3=c_2-c_1.
\end{eqnarray*}

In order to obtain not only the optimal convergence order of the error $H^1$-norm but also that of $L^2$-norm,  we uniformly refine the mesh size $h$ from $\frac{1}{4}$ to $\frac{1}{32}$ and make the time step size $\Delta t=h^2$. We also introduced an error notation $e_{u_i}=u_i(t^N)-u_i^N$ with $T=1.0, ~ N=\frac{T}{\Delta t}$ to display the results conveniently.

\begin{figure}[htbp]
\centering
\subfigure[$\kappa_1=0.01$]{
\begin{minipage}[t]{0.33\linewidth}
\centering
\includegraphics[width=45mm, height=50mm]{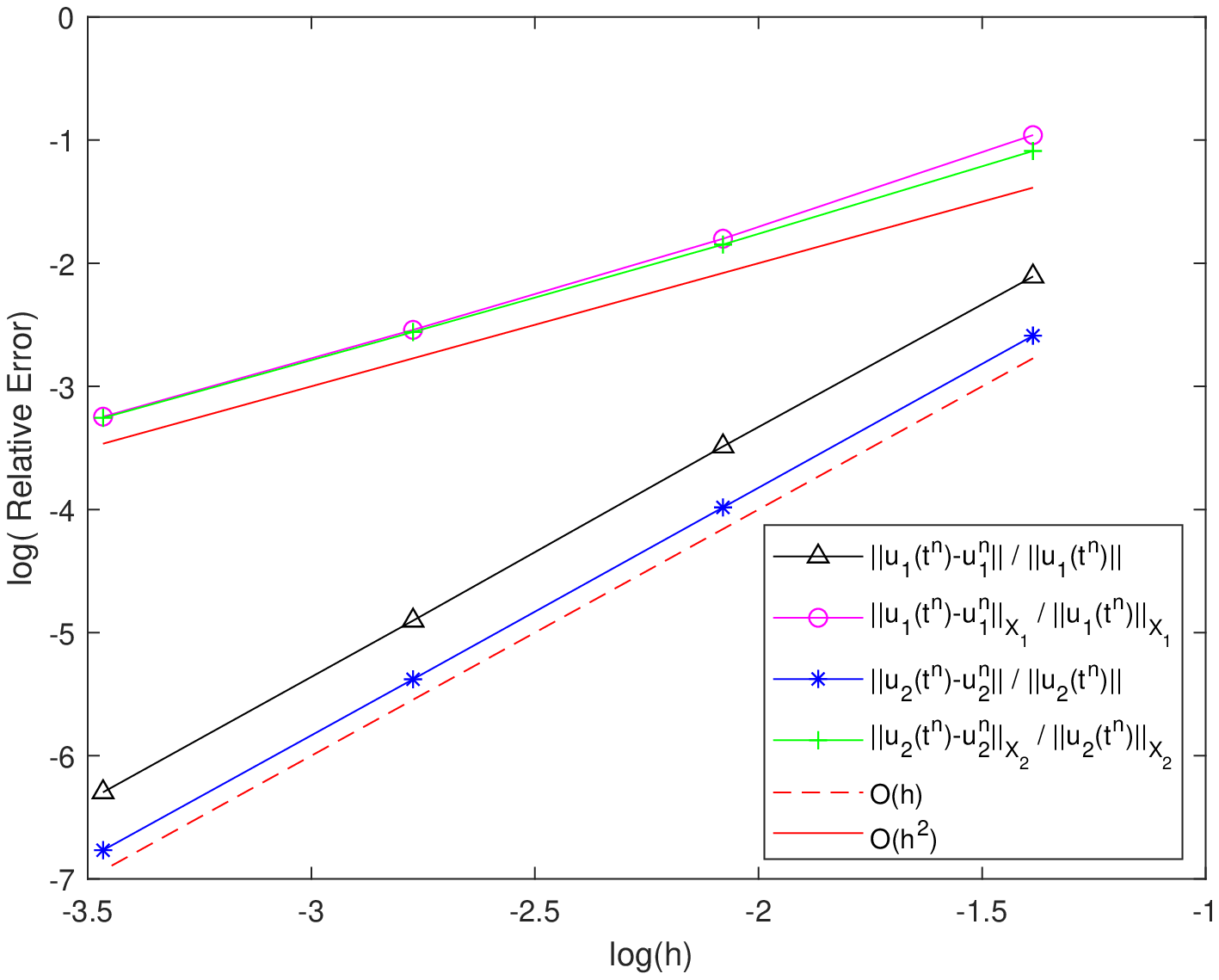}
\end{minipage}%
}%
\subfigure[$\kappa_2=1.0$]{
\begin{minipage}[t]{0.33\linewidth}
\centering
\includegraphics[width=45mm, height=50mm]{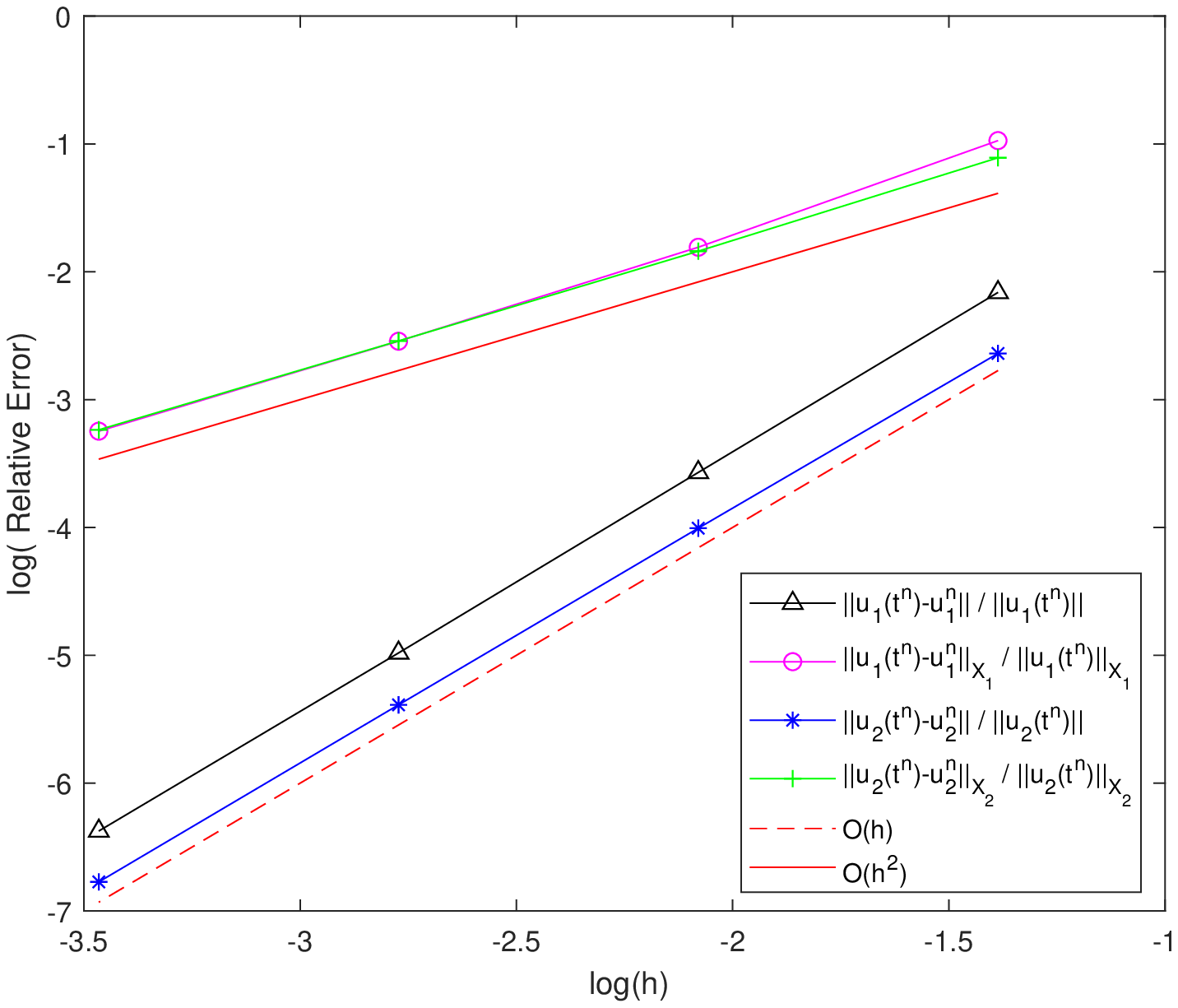}
\end{minipage}%
}%
\subfigure[$\kappa_3=10.0$]{
\begin{minipage}[t]{0.33\linewidth}
\centering
\includegraphics[width=45mm, height=50mm]{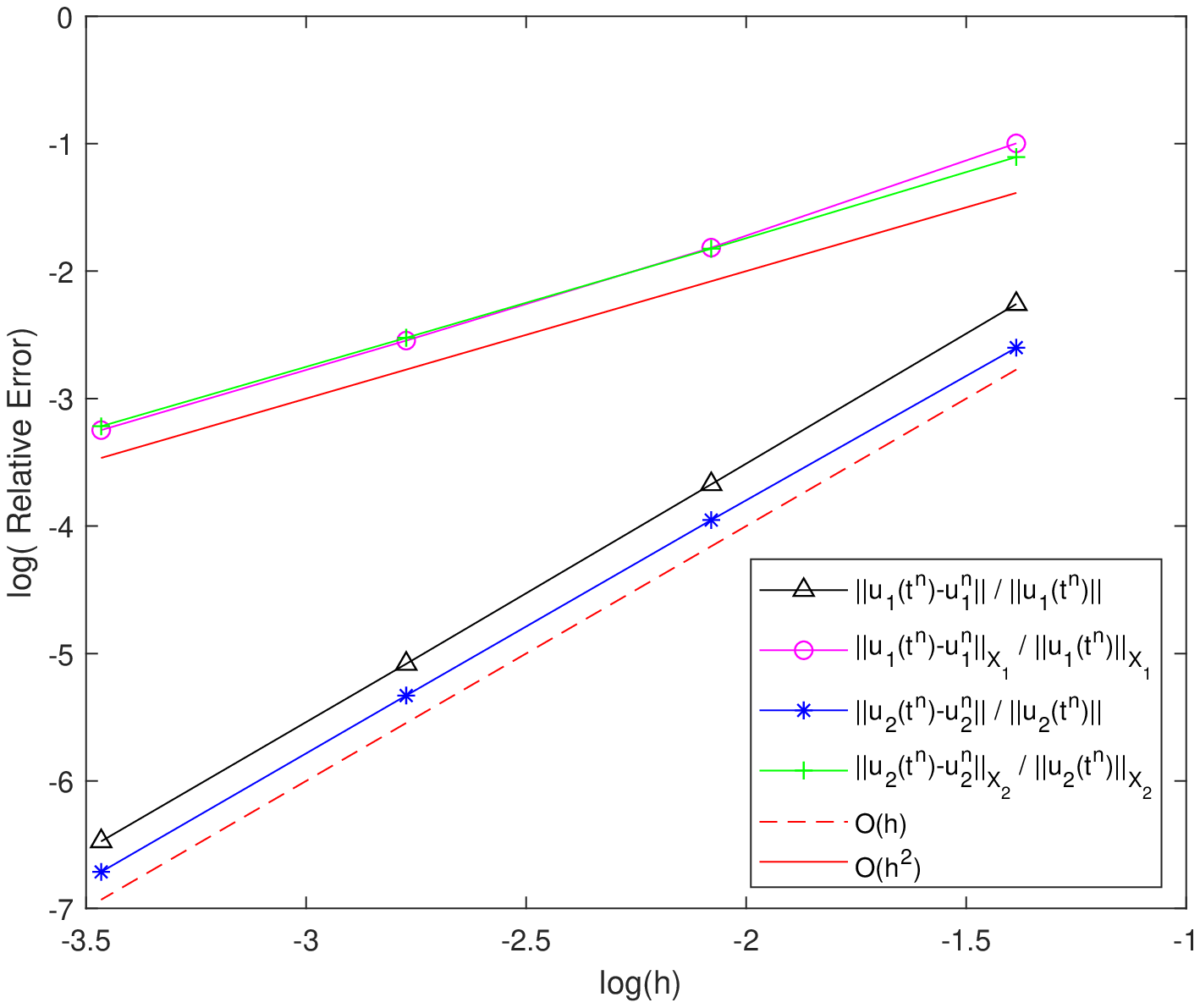}
\end{minipage}%
}%
\centering
\caption{The $L^2$-error and $H^1$-error of the $u_1$ and $u_2$ for A1 (J=3) with different $\kappa_j$.}
\label{plotk}
\end{figure}

Based on the above exact solution, we first consider the smooth heat-heat coupled problem with only a random friction parameter $\kappa$. In order to show the convergence performance better, we carry out a group of simulations with $J=3$, by which the friction parameter is selected as $\kappa_1=0.01, \kappa_2=1.0, \kappa_3=10.0$. All the other physical parameters $a,\nu_1,\nu_2$ are simply taken as $1.0$. We utilize the A1 algorithm to simulate the above CASE 1 exact solution problem, and show the convergence order of the $L^2$-norm and $H^1$-norm for two regions in Fig. \ref{plotk}. It can be observed from Fig. \ref{plotk} that the $L^2$ relative errors and $H^1$ relative errors of $u_1$ and $u_2$ for each sample have achieved the optimal convergence order. It is noteworthy that the random friction parameter $\kappa$ we choose is not a small perturbation constraint, which positively supports our advantage mentioned in Remark \ref{remk}.

Furthermore, we consider the smooth heat-heat coupled model with two random diffusion coefficients $\nu_1, \nu_2$ and one random friction parameter $\kappa$. We select a group of simulation with $J\times J, J=3$ members to test the convergence performance of A2 and A3 algorithms, including three friction parameters $\kappa_1=0.01,\kappa_2=1.0,\kappa_3=10.0$, and three diffusion coefficient groups
$\nu_{1,j}({\bf x},t)=\nu_{2,j}({\bf x},t)=1+(1+\varepsilon_j)\sin(t)$ with $\varepsilon_1=0.6207, \varepsilon_2=0.1841$ and $\varepsilon_3=0.2691$ \cite{wang zhu}. The parameter $a$ can also be taken as $1.0$. 
The approximate accuracy of the $u_1$ and $u_2$ in $L^2$-norm and $H^1$-norm of the A2 algorithm are shown in Table \ref{A21} and Table \ref{A22}, which can verify the optimal convergence orders for both $u_1$ and $u_2$. Numerical results for A3 in Table \ref{A31} and Table \ref{A32} also demonstrate very similar approximate accuracy as that of A2. Both results can indicate the effectiveness of our algorithms and confirm the theoretical analysis.
\begin{table}[!h]
\setlength{\abovecaptionskip}{0cm}
		\setlength{\belowcaptionskip}{0.2cm}
\caption{The $L^2$-error of $u_1$ and $u_2$ for A2 $(J\times J, J=3)$ with different $\kappa_j$ and $\nu_{j}$ while $\Delta t=h^2$.}
\tabcolsep 3pt \vspace*{-10pt}
\par
\begin{center}
\def\temptablewidth{1.11\textwidth}
\resizebox{\textwidth}{19mm}{
\begin{tabular*}{\temptablewidth}{@{\extracolsep{\fill}}cccccccccc}
\hline
$h$ & $\|e_{u_1}\|^{\kappa_1,\nu_1}  $& $\|e_{u_1}\|^{\kappa_1,\nu_2}$&$\|e_{u_1}\|^{\kappa_1,\nu_3}  $&$\|e_{u_1}\|^{\kappa_2,\nu_1}  $&
$\|e_{u_1}\|^{\kappa_2,\nu_2}  $&$\|e_{u_1}\|^{\kappa_2,\nu_3}  $&$\|e_{u_1}\|^{\kappa_3,\nu_1}  $&
$\|e_{u_1}\|^{\kappa_3,\nu_2} $&$\|e_{u_1}\|^{\kappa_3,\nu_3}  $\\
\hline
1/4       &0.061806      &0.053666    &0.054391     &0.059631   &0.053746     &0.053899        &0.056034    &0.055158   &0.054308   \\

1/8       &0.015217     &0.013336     &0.013488     &0.014818    &0.013350      &0.013408      &0.014019    &0.013540   &0.013416\\

1/16      &0.003791     &0.003331    &0.003367    &0.003698    & 0.003334     &0.003349        &0.003503    &0.003371  &0.003345 \\

1/32      &0.000947   &0.000832    &0.000841      &0.000924   &0.000833      &0.000837         &0.000875    &0.000842  &0.000835\\
\hline
$h$ & $\|e_{u_2}\|^{\kappa_1,\nu_1}  $& $\|e_{u_2}\|^{\kappa_1,\nu_2}$&$\|e_{u_2}\|^{\kappa_1,\nu_3}  $&$\|e_{u_2}\|^{\kappa_2,\nu_1}  $&
$\|e_{u_2}\|^{\kappa_2,\nu_2}  $&$\|e_{u_2}\|^{\kappa_2,\nu_3}  $&$\|e_{u_2}\|^{\kappa_3,\nu_1}  $&
$\|e_{u_2}\|^{\kappa_3,\nu_2} $&$\|e_{u_2}\|^{\kappa_3,\nu_3}  $\\
\hline
1/4       &0.083507     & 0.07482   &0.07679     &0.082922   &0.074048    &0.075978      &0.086434   &0.076510  &0.078567   \\

1/8       &0.020930    &0.018640  &0.019149     &0.021187      &0.018945    &0.019430       &0.022571       &0.020060   &0.020590   \\

1/16      &0.005241       & 0.004663  &0.004790   &0.005329   &0.004768    &0.004889      &0.005706     &0.005078    & 0.005211\\

1/32      & 0.001310       & 0.001166 &0.001198 &0.001334  & 0.001194   &0.001224   & 0.001430         &0.001273   &0.001306  \\
\hline
\end{tabular*}}
\end{center}\label{A21}
\end{table}

\begin{table}[!h]
\setlength{\abovecaptionskip}{0cm}
		\setlength{\belowcaptionskip}{0.2cm}
\caption{The $H^1$-error of $u_1$ and $u_2$ for A2 $(J\times J,J=3)$ with different $\kappa_j$ and $\nu_{j}$ while $\Delta t=h^2$.}
\tabcolsep 3pt \vspace*{-10pt}
\par
\begin{center}
\def\temptablewidth{1.11\textwidth}
\resizebox{\textwidth}{19mm}{
\begin{tabular*}{\temptablewidth}{@{\extracolsep{\fill}}cccccccccc}
\hline
$h$ & $\|e_{u_1}\|_{X_1}^{\kappa_1,\nu_1}  $& $\|e_{u_1}\|_{X_1}^{\kappa_1,\nu_2}$&$\|e_{u_1}\|_{X_1}^{\kappa_1,\nu_3}  $&$\|e_{u_1}\|_{X_1}^{\kappa_2,\nu_1}  $&
$\|e_{u_1}\|_{X_1}^{\kappa_2,\nu_2}  $&$\|e_{u_1}\|_{X_1}^{\kappa_2,\nu_3}  $&$\|e_{u_1}\|_{X_1}^{\kappa_3,\nu_1}  $&
$\|e_{u_1}\|_{X_1}^{\kappa_3,\nu_2} $&$\|e_{u_1}\|_{X_1}^{\kappa_3,\nu_3}  $\\
\hline
1/4       &0.310719     &0.318569    &0.316509       &0.312009    &0.320394   &0.318240       &0.315402     &0.324184   &0.322016  \\

1/8       &0.155245      &0.156257    &0.155996    & 0.155376     &0.156414        &0.156150        &0.155754        &0.156782   &0.156529\\

1/16      &0.077550     &0.077678     &0.077645      &0.077566     &0.077695     &0.077663       &0.077612      &0.077738   &0.077707  \\

1/32      &0.038764     &0.038781      &0.038776     &0.038766      &0.038782     &0.038778        &0.038772     &0.038788   &0.038784 \\
\hline
$h$ & $\|e_{u_2}\|_{X_2}^{\kappa_1,\nu_1}  $& $\|e_{u_2}\|_{X_2}^{\kappa_1,\nu_2}$&$\|e_{u_2}\|_{X_2}^{\kappa_1,\nu_3}  $&$\|e_{u_2}\|_{X_2}^{\kappa_2,\nu_1}  $&
$\|e_{u_2}\|_{X_2}^{\kappa_2,\nu_2}  $&$\|e_{u_2}\|_{X_2}^{\kappa_2,\nu_3}  $&$\|e_{u_2}\|_{X_2}^{\kappa_3,\nu_1}  $&
$\|e_{u_2}\|_{X_2}^{\kappa_3,\nu_2} $&$\|e_{u_2}\|_{X_2}^{\kappa_3,\nu_3}  $\\
\hline
1/4       &0.305074     &0.306588    &0.306167    &0.308332   &0.310523  &0.309953      &0.314882    &0.317150 &0.316561\\

1/8       &0.153796   &0.154052   &0.153984    &0.155343   & 0.155799  &0.155690     &0.158884     &0.159270 &0.159178\\

1/16      & 0.077033    &0.077081    &0.077069   &0.077811  &0.077959    &0.077926      &0.079638    &0.079755  &0.079730\\

1/32      &0.038533    &0.038546  &0.038543  &0.038922   &0.038987   &0.038973      &0.039844   &0.039893  &0.039883\\
\hline
\end{tabular*}}
\end{center}\label{A22}
\end{table}

\begin{table}[!h]
\setlength{\abovecaptionskip}{0cm}
		\setlength{\belowcaptionskip}{0.2cm}
\caption{The $L^2$-error of $u_1$ and $u_2$ for A3 $(J\times J, J=3)$ with different $\kappa_j$ and $\nu_{j}$ while $\Delta t=h^2$.}
\tabcolsep 3pt \vspace*{-10pt}
\par
\begin{center}
\def\temptablewidth{1.11\textwidth}
\resizebox{\textwidth}{19mm}{
\begin{tabular*}{\temptablewidth}{@{\extracolsep{\fill}}cccccccccc}
\hline
$h$ & $\|e_{u_1}\|^{\kappa_1,\nu_1}  $& $\|e_{u_1}\|^{\kappa_1,\nu_2}$&$\|e_{u_1}\|^{\kappa_1,\nu_3}  $&$\|e_{u_1}\|^{\kappa_2,\nu_1}  $&
$\|e_{u_1}\|^{\kappa_2,\nu_2}  $&$\|e_{u_1}\|^{\kappa_2,\nu_3}  $&$\|e_{u_1}\|^{\kappa_3,\nu_1}  $&
$\|e_{u_1}\|^{\kappa_3,\nu_2} $&$\|e_{u_1}\|^{\kappa_3,\nu_3}  $\\
\hline

1/4       &0.069406     &0.057037   &0.059113     &0.066643    &0.055908    &0.057517       &0.061175   &0.054869   &0.055456   \\

1/8       &0.017386     &0.014361   &0.014882     &0.016817   &0.014121    &0.014556    &0.015516  & 0.013758  &0.013989   \\

1/16      &0.004348      & 0.003597  &0.003726   &0.004211  &0.003539     &0.003649       & 0.003888    &0.003443   &0.003504    \\

1/32      & 0.001087    &0.000899  & 0.000932   & 0.001053   &0.000885    &0.000912       &0.000972  &0.000861  & 0.000876  \\
\hline
$h$ & $\|e_{u_2}\|^{\kappa_1,\nu_1}  $& $\|e_{u_2}\|^{\kappa_1,\nu_2}$&$\|e_{u_2}\|^{\kappa_1,\nu_3}  $&$\|e_{u_2}\|^{\kappa_2,\nu_1}  $&
$\|e_{u_2}\|^{\kappa_2,\nu_2}  $&$\|e_{u_2}\|^{\kappa_2,\nu_3}  $&$\|e_{u_2}\|^{\kappa_3,\nu_1}  $&
$\|e_{u_2}\|^{\kappa_3,\nu_2} $&$\|e_{u_2}\|^{\kappa_3,\nu_3}  $\\
\hline

1/4       &0.088783   &0.080920  &0.082745    &0.089384     & 0.081230   &0.083063       &0.094984     & 0.085424  &0.087489      \\

1/8       &0.022521  &0.020487  &0.020949    &0.023057   & 0.021069   &0.021512   &0.025003   &0.022698 &0.023201  \\

1/16      & 0.005655    &0.005145  &0.005261     &0.005812    &0.005319    &0.005428      &0.006332    &0.005761     &0.005886   \\

1/32      &0.001415   &0.001287      &0.001316   &0.001456  &0.001333     &0.001360      &0.001588  &0.001445  &0.001476 \\
\hline
\end{tabular*}}
\end{center}\label{A31}
\end{table}

\begin{table}[!h]
\setlength{\abovecaptionskip}{0cm}
		\setlength{\belowcaptionskip}{0.2cm}
\caption{The $H^1$-error of $u_1$ $u_2$ for A3 $(J\times J, J=3)$ with different $\kappa_j$ and $\nu_{j}$ while $\Delta t=h^2$.}
\tabcolsep 3pt \vspace*{-10pt}
\par
\begin{center}
\def\temptablewidth{1.11\textwidth}
\resizebox{\textwidth}{19mm}{
\begin{tabular*}{\temptablewidth}{@{\extracolsep{\fill}}cccccccccc}
\hline
$h$ & $\|e_{u_1}\|_{X_1}^{\kappa_1,\nu_1}  $& $\|e_{u_1}\|_{X_1}^{\kappa_1,\nu_2}$&$\|e_{u_1}\|_{X_1}^{\kappa_1,\nu_3}  $&$\|e_{u_1}\|_{X_1}^{\kappa_2,\nu_1}  $&
$\|e_{u_1}\|_{X_1}^{\kappa_2,\nu_2}  $&$\|e_{u_1}\|_{X_1}^{\kappa_2,\nu_3}  $&$\|e_{u_1}\|_{X_1}^{\kappa_3,\nu_1}  $&
$\|e_{u_1}\|_{X_1}^{\kappa_3,\nu_2} $&$\|e_{u_1}\|_{X_1}^{\kappa_3,\nu_3}  $\\
\hline
1/4       &0.308623    & 0.314787    &0.313126    &0.309652   &0.316327   & 0.314571       &0.312528     &0.319704    &0.317901 \\

1/8       &0.154925     & 0.155678   &0.155480     &0.155037    &0.155823    &0.155620         &0.155376      &0.156179    &0.155980   \\

1/16      &0.077509     &0.077603  &0.077578   &0.077522   &0.077620    &0.077595       &0.077565      &0.077663    &0.077639\\

1/32      & 0.038759    &0.038771    &0.038768    &0.038761    &0.038773     &0.038770       &0.038766     &0.038778     &0.038775 \\
\hline
$h$ & $\|e_{u_2}\|_{X_2}^{\kappa_1,\nu_1}  $& $\|e_{u_2}\|_{X_2}^{\kappa_1,\nu_2}$&$\|e_{u_2}\|_{X_2}^{\kappa_1,\nu_3}  $&$\|e_{u_2}\|_{X_2}^{\kappa_2,\nu_1}  $&
$\|e_{u_2}\|_{X_2}^{\kappa_2,\nu_2}  $&$\|e_{u_2}\|_{X_2}^{\kappa_2,\nu_3}  $&$\|e_{u_2}\|_{X_2}^{\kappa_3,\nu_1}  $&
$\|e_{u_2}\|_{X_2}^{\kappa_3,\nu_2} $&$\|e_{u_2}\|_{X_2}^{\kappa_3,\nu_3}  $\\
\hline
1/4       &0.304845    &0.305993  &0.305669 &0.308008 &0.309717   & 0.309272 &0.314534  &0.316144 &0.315717\\

1/8       & 0.153752 &0.153947   &0.153895 &0.155298   &0.155686   &0.155594&0.158851 &0.159152 & 0.159082  \\

1/16      &0.077028   &0.077067  &0.0770576 &0.077805 &0.077945  & 0.077914 &0.079634     & 0.079741   &0.079718 \\

1/32      &0.038532 &0.038544& 0.038541 &0.038922&  0.038985 &0.038971& 0.039843   &0.039892 &0.039882 \\
\hline
\end{tabular*}}
\end{center}\label{A32}
\end{table}

To better estimate the uncertainty and sensitivity of the solutions of random PDEs, we need to choose more samples, which may lead to a slow convergence rate. Therefore, herein, we compare the computational efficiency of our proposed algorithms A2 and A3 with the standard data-passing partitioned algorithm \cite{heat-linear} under the selected $J\times J, J=1,5,10,15,20$ realizations. Herein, for convenience, the two random diffusion coefficients can be selected as $\nu_{1}(\omega,{\bf x},t)=\nu_2(\omega,{\bf x},t)=1+(1+\omega)\sin(t)$ and the random friction parameter can be defined as $\kappa(\omega)=0.01+\omega$. The comparison of the elapsed CPU times is presented in Table \ref{cpu}, from which we can see that the proposed ensemble algorithms are meaningfully faster than the standard data-passing partitioned no-ensemble algorithm except for the case $J\times J = 1\times 1$. More importantly, both methods capture the same behaviors while the A2 saves about $91.8\%$ of the computation time, and the A3 saves about $92.1\%$ of the computation time with $J \times J = 20 \times 20$. Note that the elapsed CPU time of the A3 algorithm is not much less than that of A2, which is definitely true. That is because the A3 algorithm only calculates the average of $\nu_{i,j}$ in time on the basis of the A2 algorithm. Obviously, the elapsed CPU time for calculating more $N$ linear systems is not much longer than $ J \times J $. As has been noted, we can verify the superiority of our proposed algorithms by comparing the elapsed CPU time.
\begin{table}[!h]
\caption{The comparison of the elapsed CPU time while the mesh size $h=\frac{1}{32}$ .}
\tabcolsep 3pt \vspace*{-10pt}
\par
\begin{center}
\def\temptablewidth{1.0\textwidth}
{\rule{\temptablewidth}{1pt}}
\begin{tabular*}{\temptablewidth}{@{\extracolsep{\fill}}cccccc}
\hline
$J\times J$ & $1\times 1$ & $5 \times 5$ & $10 \times 10$ & $15 \times 15$ & $20 \times 20$\\
\hline
Standard Data-passing Partitioned &151.2   &3793.4   &15334.6  &34203.5  &60993.0\\

A2  &157.9    &442.4   &1371.6   &2886.5 &4972.2 \\

A3  &12.8    &301.1  &1198.5     &2692.9   &4801.8 \\

\hline   
\end{tabular*}%
\end{center}\label{cpu}
\end{table}

\subsection{Stability problem}
\begin{figure}[htbp]\label{energy}
\centering
{
\begin{minipage}[t]{0.45\linewidth}
\centering
\includegraphics[width=60mm, height=45mm]{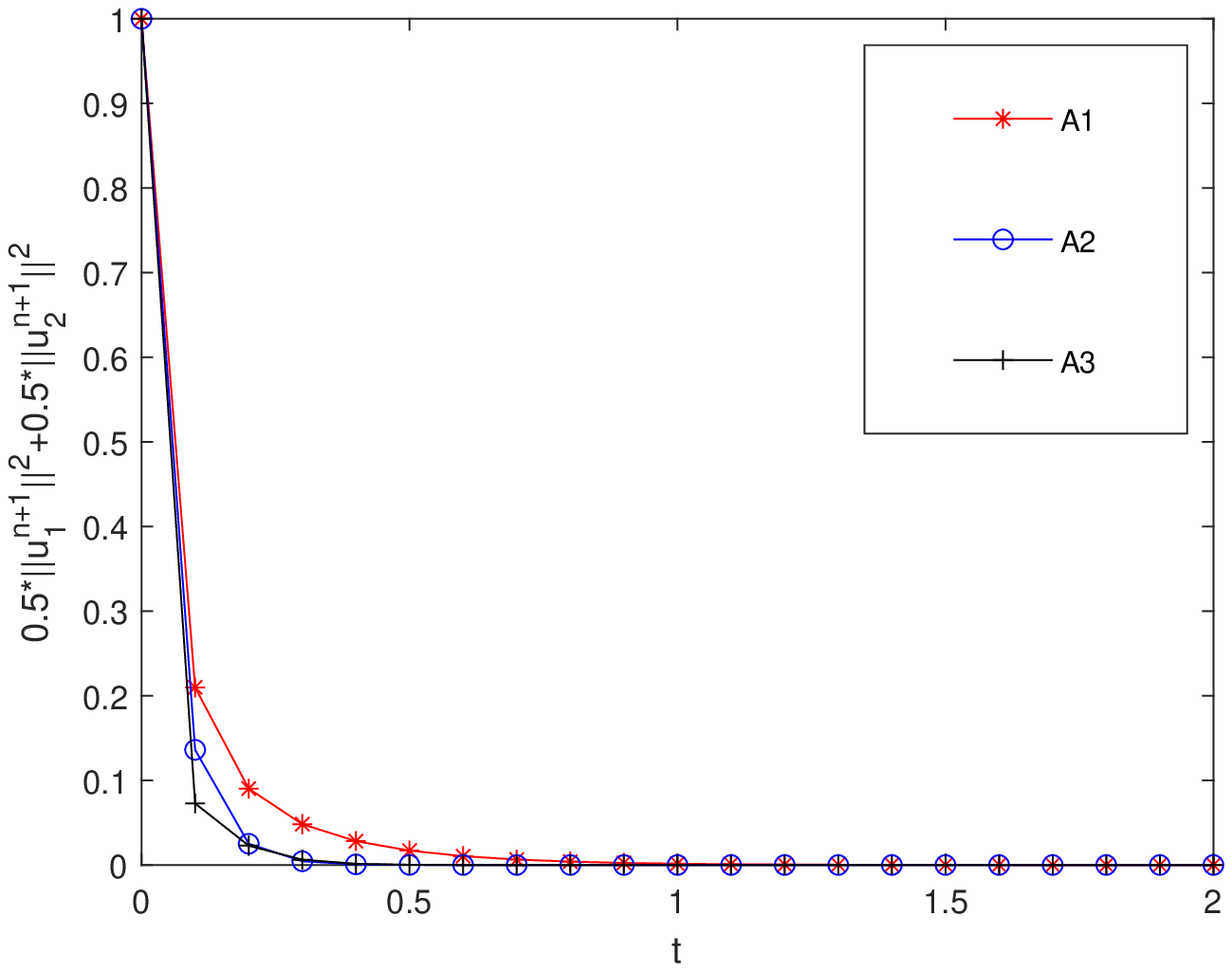}
\end{minipage}%
}%
{
\begin{minipage}[t]{0.45\linewidth}
\centering
\includegraphics[width=60mm, height=45mm]{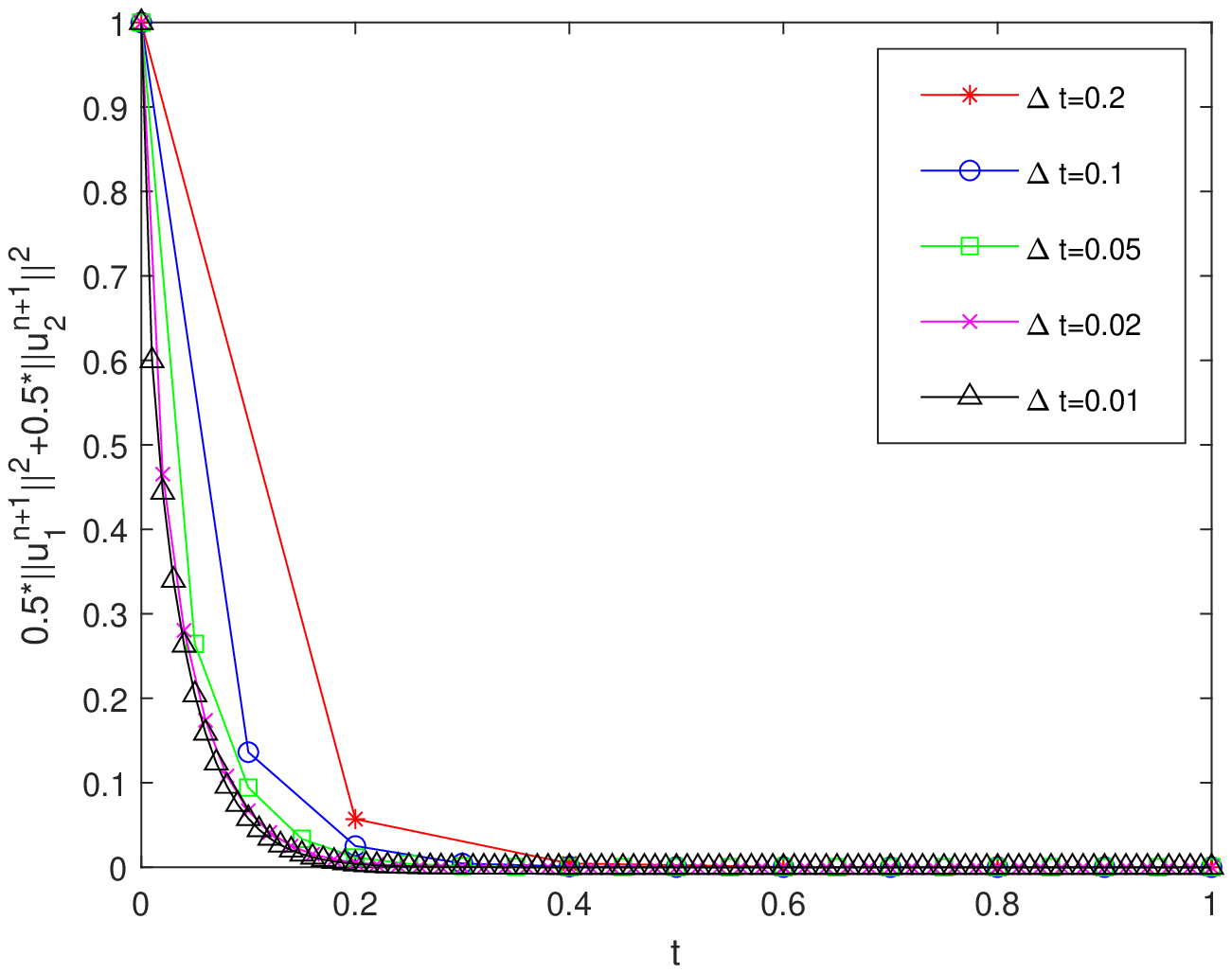}
\end{minipage}%
}%
\centering
\caption{Energy versus time of three proposed ensemble algorithms on large time step points with $\Delta t=0.1$ while $h=\frac{1}{32}$ (left), and the evolution of $0.5*\|u_i^{n+1}\|^2+0.5*\|u_2^{n+1}\|^2$ in time  for the varying time step size $\Delta t=0.2,0.1,0.05,0.02,0.01$ with fixed $h=1/32$ for the A2 algorithm (right).}
\end{figure}

In the second numerical experiment, we will test the stability of our proposed ensemble algorithms, especially unconditional stability. For simplicity, we can set $f_{i,j}=0$, $u_{i,j}^0=1.0$, and $u_{i,j}=0~\mathrm{on}~\Gamma_i$ based on the domain and parameter selections $\nu_{1}(\omega,{\bf x},t), ~ \nu_2(\omega,{\bf x},t)$ and $\kappa(\omega)$ of the first numerical experiment. Fig. \ref{energy} (left) displays the quantity of energy $0.5*{\|u_1^{n+1}\|^2}+0.5*{\|u_2^{n+1}\|^2}$ on large time step points $\Delta t=0.1$ while $h=\frac{1}{32}$, where $u_i^{n+1}=\frac{1}{J \times J} \sum_{j=1}^{J \times J}u_{i,j}^{n+1}$. We can easily calculate that the energy equals $1$ while $t=0$. Then, it is clear to see that the energy of both algorithms decreases rapidly and finally reaches a stable state.

Moreover, to better display that the proposed ensemble algorithms are unconditionally stable, we check the energy attenuation results of A2 as a representative under different time steps. As illustrated in Fig. \ref{energy} (right), we plot the evolution of $0.5*{\|u_1^{n+1}\|^2}+0.5*{\|u_2^{n+1}\|^2}$ in time for the varying time step size $\Delta t=0.2,0.1,0.05,0.02,0.01$ with fixed mesh size $\frac{1}{32}$. We note that, for fixed $h$, with the different time steps $\Delta t$, the energy of A2 decreases rapidly with the increase of time and tends to be stable. More importantly, the smaller the time steps $\Delta t$, the quicker we reach the steady case. Generally, the observation shows that the time step conditions are not necessary, which can indicate that the proposed algorithms are unconditionally stable.

\subsection{Steel-titanium composite plate fuel cell model}
Due to the composite materials having environmental protection, low cost and excellent mechanical properties, it's have been widely used in life and production, such as steel-titanium composite plates, copper-aluminium composite plates, etc. It is well-known that such composite plates are the critical components of the metal bipolar plate material for the fuel cell. The performance and life of the cell will be affected by temperature, and the thermal management of the cell has received extensive attention. Therefore, based on the above application, we design a steel-titanium composite plate fuel cell model to observe its internal heat conduction phenomenon. 

\begin{figure}[htbp]\label{3dmodel}
\centering
\includegraphics[width=70mm,height=40mm]{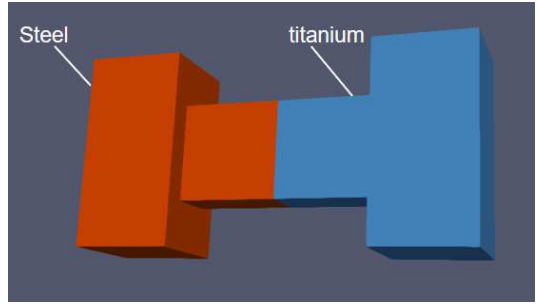}
\caption{The conceptual domain of steel-titanium composite plate full cell model.}
\end{figure}
We assume two materials can be expressed explicitly as $[0,2]\times[0,2]\times[0,4]\cup[2,4]\times[0.5,1]\times[2,3]\cup[6,8]\times[0,2]\times[2,4]\cup[4,6]\times[0.5,1]\times[2,3]$ in the conceptual domain, see Fig. \ref{3dmodel}. For simplicity, we only give the leftmost heat of the steel plate as Dirichlet boundary condition $u_{1}(t)=20$ and assume other surfaces of the whole model (excluding the interface)  as homogeneous Neumann boundary conditions $\nabla u_{i}(t) \cdot \hat{n}_i=0$. The initial conditions are presumed to be $0$ and there is no external force term. Owing to the random fiction parameter $\kappa$  does not need a small perturbation constraint, we select ten samples for the random $\kappa$, which are $0.001, 0.005, 0.01, 0.05, 0.1, 0.5, 1, 5, 10, 50$. According to Wikipedia, the heat conductivity of steel is usually between $51.9$ and $67.4$ due to the different carbon content, and the heat conductivity of titanium is about $15.6-22.5$. Hence, in the above selection interval, we randomly selected ten samples for the Monte Carlo simulation.

\begin{figure}[htbp]
\centering
\subfigure[$t=0.05$]{
\begin{minipage}[t]{0.5\linewidth}
\centering
\includegraphics[width=65mm, height=40mm]{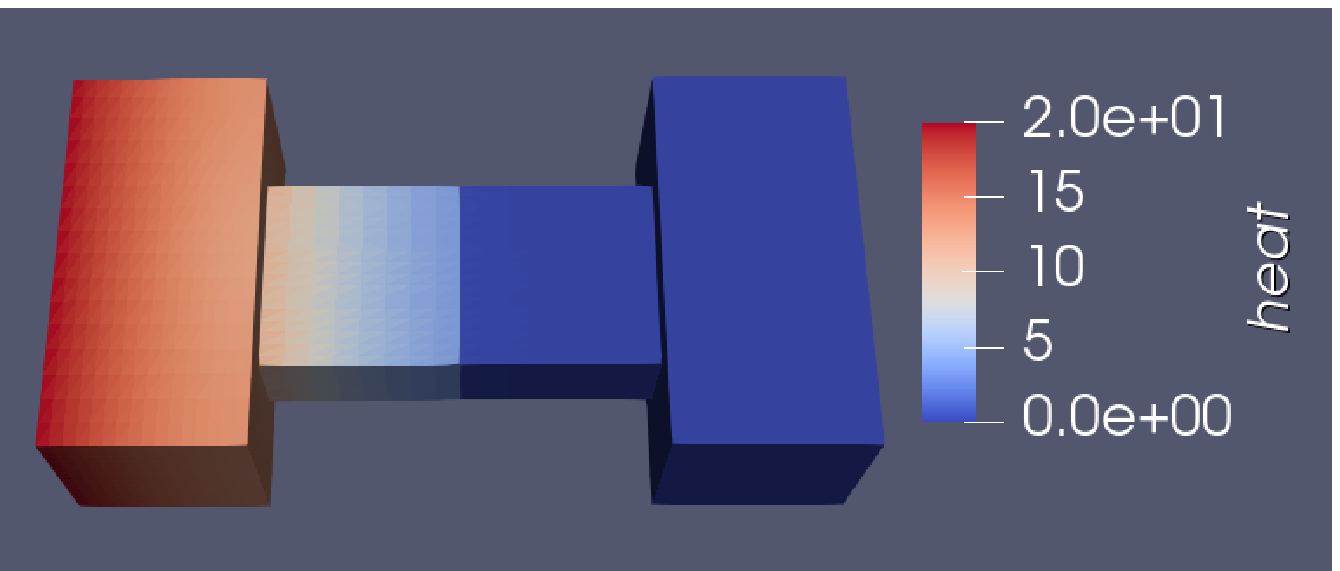}
\end{minipage}%
}%
\subfigure[$t=0.2$]{
\begin{minipage}[t]{0.5\linewidth}
\centering
\includegraphics[width=65mm, height=40mm]{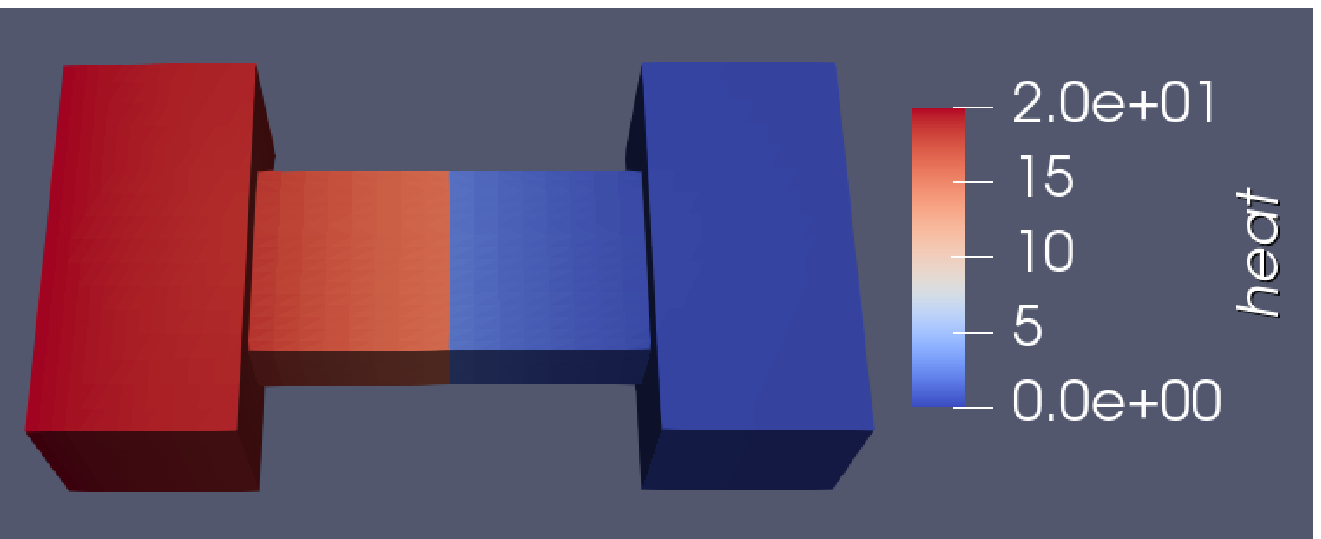}
\end{minipage}%
}%
\\
\subfigure[$t=1.0$]{
\begin{minipage}[t]{0.5\linewidth}
\centering
\includegraphics[width=65mm, height=40mm]{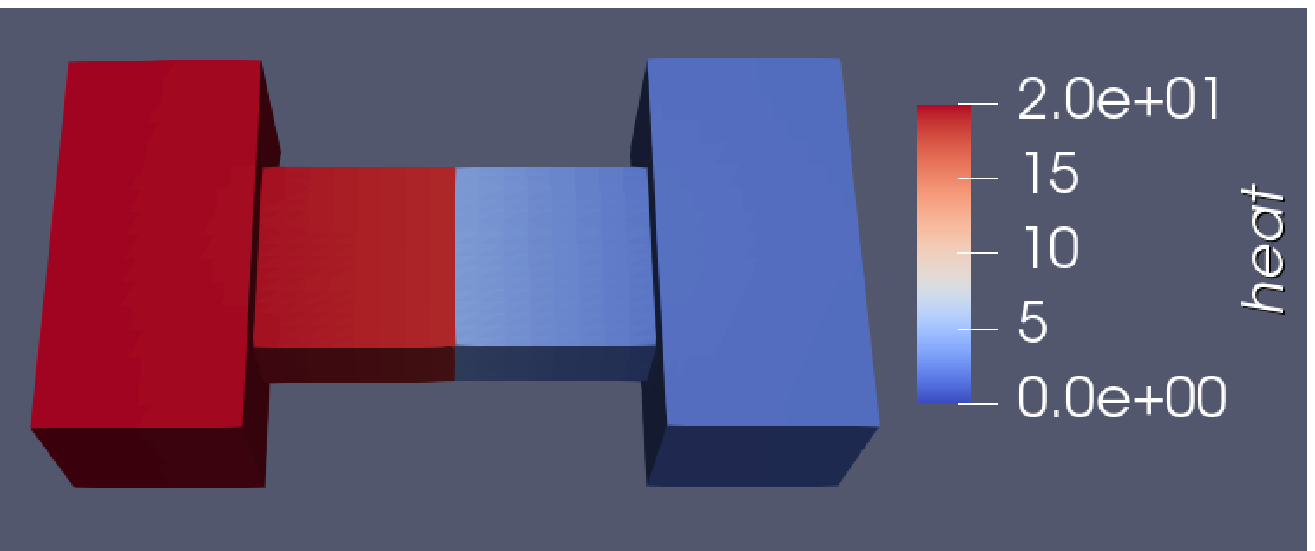}
\end{minipage}%
}%
\subfigure[$t=4.0$]{
\begin{minipage}[t]{0.5\linewidth}
\centering
\includegraphics[width=65mm, height=40mm]{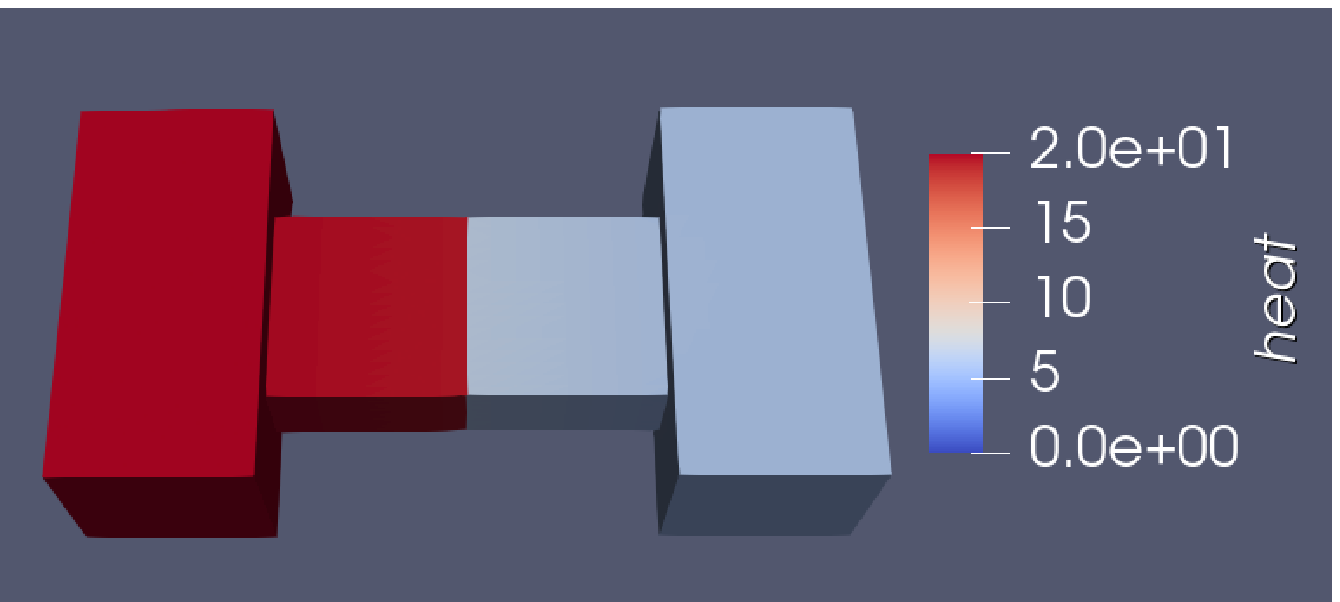}
\end{minipage}%
}%
\centering
\caption{Heat conduction in the steel-titanium composite plate fuel cell model with $h=\frac{1}{16}$ and $\Delta t=0.01$.}
\label{heat}
\end{figure}

We utilize the A2 algorithm with $h=\frac{1}{16}$ and $\Delta t=0.01$ to simulate the heat transfer process inside the steel-titanium composite plate fuel cell model, as depicted in Fig. \ref{heat}. With the increase of time, the heat is gradually transferred from the steel plate to the titanium plate, and because the steel plate has a relatively high heat conductivity, it can be heated rapidly. Moreover, there is a significant jump at the interface between the steel plate and the titanium plate influenced by the friction parameter. In summary, all discussions further guarantee the effectiveness of our proposed algorithms.

\section{Conclusion}
In this paper, three efficient ensemble algorithms are proposed to solve the heat-heat coupled system with two random diffusion coefficients $\nu_1, \nu_2$, and one random friction parameter $\kappa$. 
The main idea of this paper is to combine the Monte Carlo method with the ensemble idea to solve the random parabolic PDEs. For simplicity, we first consider the heat-heat problem with the random friction parameter $\kappa$ only and present an unconditionally stable and convergent A1 algorithm. Then, an advanced heat-heat coupled problem with three random coefficients is studied. We construct two ensemble algorithms for fast solving this more complex system and confirm that they are unconditionally stable and have optimal convergence order under a parameter condition. This kind of idea can also be extended to the nonlinear fluid-fluid problem.

\end{document}